\documentclass[11pt,review]{article}
\usepackage{rotating}
\usepackage{amssymb}
\usepackage{color}
\usepackage{amsthm}
\usepackage{amsfonts}
\usepackage{graphicx}
\usepackage{graphics}
\usepackage{setspace}
\usepackage{bm}
\usepackage{amsmath}
\usepackage{algorithm}
\usepackage{algpseudocode}
\usepackage{url}
\usepackage{url}
\usepackage[format=plain,justification=raggedright,singlelinecheck=false]{caption}
\usepackage[justification=centering]{caption}
\makeatletter
\def\BState{\State\hskip-\ALG@thistlm}
\makeatother \numberwithin{equation}{section}

\newenvironment{Proof}{\noindent{\bf {Proof}.}}{\hfill{$\blacksquare$}}



\newcommand{\e}{\mathrm{e}}                         %

\usepackage{color}
\usepackage[draft,disable]{todonotes}
\newcommand{\TODO}[1]{\textcolor{black}{#1}}
\newcommand{\DO}[1]{\textcolor{black}{#1}}
\newcommand{\TO}[1]{\textcolor{black}{#1}}
\begin{document}
\date{\small\textsl{\today}}
\title{Generalized moving least squares and moving
 kriging least squares approximations for solving the transport equation on the
sphere}
\author{\large
Vahid Mohammadi  $^{\mbox{\footnotesize a}}$, \large  Mehdi Dehghan
$^{\mbox{\footnotesize a}}$, \large Amirreza Khodadadian $^{\mbox{\footnotesize b}}$\footnote{Corresponding author,~ {\em E-mail addresses:}
	v.mohammadi@aut.ac.ir, v.mohammadi.aut@gmail.com (V. Mohammadi),
	mdehghan@aut.ac.ir, mdehghan.aut@gmail.com (M. Dehghan),
	amirreza.khodadadian@ifam.uni-hannover.de (A. khodadadian),
	thomas.wick@ifam.uni-hannover.de (T. Wick)} 
, \large Thomas Wick $^{\mbox{\footnotesize b}}$
\\ \\
\small{\em $^{\mbox{\footnotesize a}}$\em Department of Applied
Mathematics, Faculty of Mathematics
and Computer Sciences,}\vspace{-1mm}\\
\small{\em Amirkabir University of Technology, No. 424, Hafez
Ave.,15914, Tehran, Iran}\\ \\
\small{\em $^{\mbox{\footnotesize b}}$\em Institute of Applied Mathematics, Leibniz University of Hannover,
}\vspace{-1mm}\\
\small{\em  Welfengarten 1, 30167
    Hanover, Germany}}
\maketitle
 \vspace{.1cm}
\bigskip
\leftskip=0pt \hrule\vskip 8pt \hspace*{-.8cm} {\bfseries Abstract.}
In this work, we apply two meshless methods for the numerical
solution of the time-dependent transport equation defined on the
sphere in spherical coordinates. The first technique, which was
introduced by Mirzaei (BIT Numerical Mathematics, 54 (4) 1041-1063, 2017) in Cartesian coordinates is a
generalized moving least squares approximation, and the second one,
which is developed here, is moving kriging least squares
interpolation on the sphere. These methods do not depend on the
background mesh or triangulation, and they can be implemented on the
transport equation in spherical coordinates easily using different
distribution points. Furthermore, the time variable is approximated
by a second-order backward differential formula. The obtained fully
discrete scheme is solved via the biconjugate gradient stabilized
algorithm with zero-fill incomplete lower-upper (ILU) preconditioner
at each time step. Three well-known test problems namely solid body
rotation, vortex roll-up, and deformational flow are solved
to demonstrate our developments.

\bigskip
\leftskip=0pt \hrule\vskip 8pt \hspace*{-.8cm}
 \hspace*{.4cm}
\hspace*{.5cm}

\section{Introduction}
 In this paper, we consider the following
transport equation that is defined on the sphere in spherical
coordinates \cite{flyer2007transport,fornberg2011stabilization}
\begin{equation}\label{Eq-1}
\dfrac{\partial u}{\partial t}+\bm{v}.\nabla u=0,
\end{equation}
where $u$ is the scalar quantity being transported \TODO{and
 $\bm{v}=\bm{v}(\lambda,\theta,t)=(v_{1}(\lambda,\theta,t),v_{2}(\lambda,\theta,t))^{T}$
 represents the velocity field. Here $ \lambda$ denotes the longitude and $\theta$ is the latitude
which both are measured from the equator \cite{flyer2007transport}.} Furthermore,
 $\nabla$ is the gradient operator on the surface of a unit
sphere in spherical coordinates which is defined as
\cite{flyer2007transport,fornberg2011stabilization}
\begin{equation}\label{Eq-2}
\nabla:=\left(\dfrac{1}{\cos(\theta)} \dfrac{\partial}{\partial
\lambda},  \dfrac{\partial}{\partial \theta}\right)^{T}.
\end{equation}
As seen, in the north and south pole, i.e., $\theta=\pm
\dfrac{\pi}{2}$, this operator is singular.

\TODO{The transport equation has various applications such as modeling transport in layered magnetic
    materials \cite{maclaren2000first}, spin valves
    \cite{krems2007boltzmann}, ocean surface modeling
    \cite{bender1996modification}, numerical weather prediction
    \cite{cotter2012mixed}, modeling oil weathering, and transport in sea
    ice \cite{afenyo2016modeling}. Moreover, we should note that in the atmospheric modeling, transport processes have significant
    importance \cite{nair2010class}.}

  Different standard tests are
considered for this problem, which are solid body rotation and
deformational flow \cite{nair2010class}. A benchmark transport test
on the sphere is the solid body rotation of a cosine bell along a
great circle trajectory was introduced in
\cite{williamson1992standard}. Another test, which is introduced and
studied in \cite{nair2008moving} is deformational flow (vortex). In
the Cartesian geometry, two important deformational tests with the
analytic solutions are called the ''Smolarkiewicz's test''
\cite{smolarkiewicz1982multi} and ''Doswell vortex''
\cite{doswell1984kinematic}, respectively (see also
\cite{staniforth1987comments}). Besides, Leveque introduced a
deformational test, in which flow trajectories are much more complex
\cite{leveque1996high}. In \cite{nair2010class}, the authors
followed \cite{leveque1996high}, and constructed different
deformational tests in Cartesian and spherical geometries, in which
some tests are non-divergent flow and the others are divergent. As
said in \cite{flyer2007transport}, geophysical fluid motions on all
scales are dominated by the advection process. Therefore, computing
the numerical solutions play a more important role for solving the
transport equation.

In recent years, there are diverse research works based on numerical
methods to solve the transport equation on the sphere such as
high-order finite volume methods (FVMs)
\cite{cheruvu2007spectral,chen2008shallow,zerroukat2004slice},
continuous and discontinuous Galerkin (DG) methods
\cite{giraldo2000lagrange,nair2005discontinuous,taylor2007mass},
radial basis functions (RBFs) in spherical coordinates
\cite{flyer2007transport}, radial basis functions
\cite{flyer2009radial}, adaptive mesh refinement technique
\cite{jablonowski2006block, lauter2007parallel,st2008comparison},
\TO{the conservative semi-Lagrangian multi-tracer transport scheme
\cite{lauritzen2012standard},} stabilization of RBF-generated finite
difference techniques (RBF-FD) in spherical coordinates (adding an
artificial hyperviscosity) \cite{fornberg2011stabilization} as well
as global, local and partition of unity RBFs methods combined with
the semi-Lagrangian approach in Cartesian coordinates
\cite{shankar2018mesh},
 and a
higher-order compatible finite element scheme for the nonlinear
rotating shallow water equations on the sphere
\cite{shipton2018higher}.

In this article, we employ two meshless techniques, namely
generalized moving least squares (GMLS) and moving kriging least
squares (MKLS) approximations on the sphere in spherical coordinates
for discretizing the spatial variables of the transport equation
(\ref{Eq-1}). For the first time, the GMLS technique in subdomains
of $\mathbb{R}^d$ was introduced by Mirzaei and his co-workers
\cite{mirzaei2012generalized}. The generalized moving least squares
reproducing kernel approach is also introduced and analyzed in
$\Omega \subset \mathbb{R}^d$ \cite{salehi2013generalized}.
Recently, the GMLS technique on the sphere was introduced and
analyzed by Mirzaei in Cartesian coordinates
\cite{mirzaei2017direct}. Here, we approximate $\nabla u$ defined in
Eq. (\ref{Eq-1}) via GMLS in spherical coordinates. Besides, we have
developed an MKLS interpolation on the sphere, which approximates
$\nabla u$ in spherical coordinates. This technique was first
introduced in \cite{gu2003moving} for subdomains in $\mathbb{R}^d$,
and it is not considered on the sphere. Against GMLS approximation,
the MKLS method satisfies the Kronecker the delta property
\cite{gu2003moving}. Of course, it depends on a parameter that is
similar to the shape parameter of an RBF interpolation. As discussed
in \cite{mirzaei2017direct,schaback2017error}, these methods can be
considered as "generalized finite differences" in which the
differential operators involved in a PDE such as Eq. (\ref{Eq-1}),
can be approximated at each scattered data point on each local
sub-domain. Both presented approximations are implemented simply on
transport equation defined on the sphere since they do not depend on
a background mesh or triangulation.

\TODO{The main advantage of
GMLS and MKLS approximations developed here is that
since two techniques do not depend on any background
mesh or triangulation, there is no difficulty in implementation. In this work, we apply them on a transport equation on
the unit sphere via two different set of points.
The proposed methods can be simply implemented to solve numerically various model equations defined on the sphere in different scientific problems.}
The temporal variable of Eq. (\ref{Eq-1}) is discretized by a
second-order backward differential formula (BDF)
\cite{li2018second,li2017second}.

The remainder of this manuscript
is as follows. In Section \ref{Sec-2}, the time variable of Eq.
(\ref{Eq-1}) is discretized by a second-order backward differential
formula. In Section \ref{Sec-3}, the generalized moving least
squares approximation in spherical coordinates, and how it can be
applied to approximate the advection operator in Eq. (\ref{Eq-1})
are presented. In Section \ref{Sec-4}, a new approximation namely
moving kriging least squares is introduced on the sphere, and we
have approximated the advection operator of an unknown solution
using this technique. In Section \ref{Sec-5}, we have obtained the
full-discrete scheme of transport equation (\ref{Eq-1}) defined on
the unit sphere due to the time and spatial discretizations proposed
here. Some numerical simulations are reported in Section \ref{Sec-6}
for three test problems, which were studied in the literature works.
Finally,
 concluding remarks are given
in Section \ref{Sec-7}.

\section{The time discretization}\label{Sec-2}
In this section, we apply a second-order
BDF for discretization the time variable of Eq. (\ref{Eq-1})
\cite{li2018second,li2017second}.
For this purpose, the time interval $[0,T]$ is divided uniformly into $M$
sub-intervals  such that $T=M\Delta t$, where $\Delta t$ indicates the time
step. By defining  $t_n:=n\Delta t$ and $u^{n}:=u(t_n)$, a second-order BDF for transport equation (\ref{Eq-1}) can be written
as follows
\begin{equation}\label{BDF-1}
\dfrac{{3{u^{n + 1}} - 4{u^n} + {u^{n - 1}}}}{{2\Delta t}} +
\dfrac{{v_1^{n + 1}}}{{\cos (\theta )}}\dfrac{{\partial {u^{n +
1}}}}{{\partial \lambda }} + v_2^{n + 1} \dfrac{{\partial {u^{n +
1}}}}{{\partial \theta }} = 0,\,\,\,\,\,\,n = 1,2,...,M - 1,
\end{equation}
\TODO{where $v_1^{n + 1}$ and $v_2^{n + 1}$ are the components of
the velocity field at $t=t_{n+1}$.}

 Also, for the first time step, we have used
the first step of backward time stepping as follows
\begin{equation}\label{BDF-2}
\dfrac{{{u^1} - {u^0}}}{{\Delta t}} + \dfrac{{v_1^1}}{{\cos (\theta
)}}\dfrac{{\partial {u^1}}}{{\partial \lambda }} +
v_2^1\frac{{\partial {u^1}}}{{\partial \theta }} = 0,
\end{equation}
\TODO{where $v_1^{1}$ and $v_2^{1}$ are the components of the
velocity field at $t=t_{1}$.}

Reformulating Eq. (\ref{BDF-1}), we have
\begin{equation}\label{BDF-3}
3{u^{n + 1}} + 2\Delta t\left( {\dfrac{{v_1^{n + 1}}}{{\cos (\theta )}}\dfrac
{{\partial {u^{n + 1}}}}{{\partial \lambda }} + v_2^{n + 1}\dfrac{{\partial {u^{n + 1}}}}
{{\partial \theta }}} \right) = 4{u^n} - {u^{n - 1}},\,\,\,\,\,\,n = 1,2,...,M - 1.
\end{equation}
Besides, Eq. (\ref{BDF-2}) can be rewritten in the following formula
\begin{equation}\label{BDF-4}
{u^1} + \Delta t\left( {\dfrac{{v_1^1}}{{\cos (\theta )}}
\dfrac{{\partial {u^1}}}{{\partial \lambda }} + v_2^1\dfrac{{\partial {u^1}}}{{\partial \theta }}} \right) = {u^0}.
\end{equation}
In what follows, we will come back to Eqs. (\ref{BDF-3}) and
(\ref{BDF-4}) for deriving their full-discrete schemes.

\section{The GMLS formulation for advection operator}\label{Sec-3}
As was mentioned earlier, the GMLS approximation on the sphere was
introduced by Mirzaei in his recent work \cite{mirzaei2017direct}.
Here, we have the derivation of \TODO{advection operator
(\ref{Eq-2}) of a given function such as $u$
 in spherical coordinates via the
GMLS approximation.}

Assume that $u \in C^{m+1}(\mathbb{S}^2)$ is a function defined on
the unit sphere $\mathbb{S}^2$, \TODO{where $\mathbb{S}^2=\{(x,y,z)
\in \mathbb{R}^3\,|\, x^2+y^2+z^2=1 \}$}. Besides, consider a set of $N$
points $X=\{\bm{x}_{1},\bm{x}_{2},...,\bm{x}_{N}\}$ on
$\mathbb{S}^{2}$. The approximation of $u$ by GMLS can be written
\cite{mirzaei2017direct}
\begin{equation}\label{GMLS-1}
u(\bm{x}) \approx \overline{u(\bm{x})}=\displaystyle \sum_{j \in I(\bm{x})} a_{j}(\bm{x})u_{j},\,\,\,\, \bm{x}
\in \mathbb{S}^{2},
\end{equation}
\TODO{where $u_{j}$ is the value $u$ at point $\bm{x}_{j} \in
I(\bm{x})$.} Here $I(\bm{x})$ is a set of indices for scattered points $X$ that
is defined as \cite{mirzaei2017direct,wendland2001moving}
$$I(\bm{x}):=\{j \in \{1,2,...,N\}: dist(\bm{x},\bm{x}_{j})< \delta\},$$
of centers contained in the cap of radius $\delta>0$ around $\bm{x}
\in \mathbb{S}^{2}$, \DO{and $dist(\bm{x},\bm{x}_{j})$ represents
geodesic distance between $\bm{x}$ and $\bm{x}_{j}$.} In spherical
coordinates, $\bm{x}:=(x,y,t)^{T} \in \mathbb{S}^2$ can be
considered as follows, e.g., see \cite{flyer2007transport}
\begin{equation}\label{GMLS-2}
x=\cos(\lambda)\cos(\theta),\,\,\,\, y=\sin(\lambda)\cos(\theta),\,\,\,\, z=\sin(\theta).
\end{equation}
In Eq. (\ref{GMLS-1}), $a_{j}(\bm{x}),\,\,j=1,2,...,|I(\bm{x})|$ for each point $\bm{x} \in \mathbb{S}^2$, is constructed
in the following vector form \cite{mirzaei2017direct}
\begin{equation}\label{GMLS-3}
\bm{a}^{\star}(\bm{x})=W(\bm{x})P(\bm{x})\Big(P^{T}(\bm{x})W(\bm{x})P(\bm{x})\Big)^{-1}Y(\bm{x})
\TODO{=:\left[a_1(\bm{x}),\cdots, a_{|I(\bm{x})|}(\bm{x})
\right]^{T}}.
\end{equation}

Here $Y$ is the vector of spherical harmonics of degree at most $m$
\cite{atkinson2012spherical,mirzaei2017direct,wendland2001moving},
\TODO{$P(\bm{x}) \in \mathbb{R}^{|I(\bm{x})| \times N(3,m)}$ is a
matrix, and its rows contain the vector of spherical harmonics.}
$W(\bm{x})$ represents a diagonal matrix with size $|I(\bm{x})|
\times |I(\bm{x})|$ with elements
$\phi\Big(\dfrac{dist(\bm{x},\bm{y})}{\delta}\Big)$
 on its diagonal, where $\phi:[0,\infty)
\rightarrow [0,\infty)$ and with the $\phi(r)>0$ for $r \in [0,1/2]$
and $\phi(r)=0$ for $r\geq 1$
\cite{mirzaei2017direct,wendland2001moving}, or
\begin{equation}\label{GMLS-4}
w(\bm{x},\bm{y}):=\phi
\Big(\dfrac{dist(\bm{x},\bm{y})}{\delta}\Big),~~~~~~ \bm{x}, \bm{y}
\in \mathbb{S}^{2},
\end{equation}
and $\delta >0$, and it is known as the weight function. As mentioned in
\cite{mirzaei2017direct,wendland2001moving}, there are different weight functions as Eq.
(\ref{GMLS-4}),
which can be considered
in this approximation. In this article, the following weight function has been used in GMLS approximation
\cite{fasshauer2007meshfree,wendland2004scattered}
\[\phi (r) = \left\{ \begin{array}{l}
{(1 - r)^4}(4r + 1),\,\,\,\,\,\,\,\,\,\,\,\,\,\,\,\,\,\,\,\,\,\,0 \le r \le 1,\\\\
0,\,\,\,\,\,\,\,\,\,\,\,\,\,\,\,\,\,\,\,\,\,\,\,\,\,\,\,\,\,\,\,\,\,\,\,\,\,\,\,\,\,
\,\,\,\,\,\,\,\,\,\,\,\,\,\,\,\,\,\,\,\,r> 1,
\end{array} \right.\]
where $r={dist(\bm{x},\bm{y})}/{\delta}$ and $\bm{x},\bm{y} \in
\mathbb{S}^2$.

If in Eq. (\ref{GMLS-1}) we consider the spherical gradient, i.e., $\nabla_{\mathbb{S}^2}:=\nabla_{0}$,
we will have \cite{mirzaei2017direct}
\begin{equation}\label{GMLS-5}
\nabla_{0}u(\bm{x}) \approx \overline{\nabla_{0}u(\bm{x})}=\displaystyle \sum_{j \in I(\bm{x})}
a_{j,\nabla_{0}}(\bm{x})u_{j},\,\,\,\, \bm{x}
\in \mathbb{S}^{2},
\end{equation}
where
\begin{equation}\label{GMLS-6}
\bm{a}^{\star}_{\nabla_{0}}(\bm{x})=W(\bm{x})P(\bm{x})\Big(P^{T}(\bm{x})W(\bm{x})P(\bm{x})\Big)^{-1}{\nabla_{0}}(Y(\bm{x})),
\end{equation}
in which $\nabla_{0}$ acts only on the vector of spherical
harmonics of degree at most $m$, i.e., $Y(\bm{x})$, $\bm{x} \in
\mathbb{S}^2$. Due to \cite[Lemma 3.1]{mirzaei2017direct} or \cite[Definition 3.3]{mirzaei2017direct}, it is easy to compute
the surface gradient in spherical coordinates. Therefore, the partial
derivatives of $a_{j}(\bm{x})$ with respect to
$\lambda$ and $\theta$ can be written as follows
\begin{align}\label{GMLS-7}\nonumber
\dfrac{{\partial {a_j}(\bm{x})}}{{\partial \lambda }} &=
{{{\left( {{\nabla _0}} \right)}_1}{a_j}(\bm{x})} \dfrac{{\partial x}}{{\partial \lambda }}
+  {{{\left( {{\nabla _0}} \right)}_2}{a_j}(\bm{x})} \dfrac{{\partial y}}{{\partial \lambda }}
+  {{{\left( {{\nabla _0}} \right)}_3}{a_j}(\bm{x})} \dfrac{{\partial z}}{{\partial \lambda }}\\\nonumber
\\
 &=  {{{\left( {{\nabla _0}} \right)}_1}{a_j}(\bm{x})}
\left( { - \sin (\lambda )\cos (\theta )} \right) + {{{\left( {{\nabla _0}} \right)}_2}{a_j}(\bm{x})}
\left( {\cos (\lambda )\cos (\theta )} \right),
\end{align}

\begin{align}\label{GMLS-8}\nonumber
\dfrac{{\partial {a_j}(\bm{x})}}{{\partial \theta }} &=
 {{{ {{\left(\nabla _0\right)}}}_1}{a_j}(\bm{x})}\dfrac{{\partial x}}{{\partial \theta }}
 + {{{\left( {{\nabla _0}} \right)}_2}{a_j}(\bm{x})} \dfrac{{\partial y}}{{\partial \theta }}
 +  {{{\left( {{\nabla _0}} \right)}_3}{a_j}(\bm{x})}\dfrac{{\partial z}}{{\partial \theta }}\\\nonumber
\\
 &=  {{{\left( {{\nabla _0}} \right)}_1}{a_j}(\bm{x})}\left(
 { - \cos (\lambda )\sin (\theta )} \right) +{{{\left( {{\nabla _0}} \right)}_2}{a_j}(\bm{x})} \left
 ( { - \sin (\lambda )\sin (\theta )} \right) +  {{{\left( {{\nabla _0}} \right)}_3}{a_j}(\bm{x})}
 ( {\cos (\theta )}),
\end{align}
where $(\nabla_{0})_{1}$, $(\nabla_{0})_{2}$ and $(\nabla_{0})_{3}$ are the components of $\nabla_{0}$. Inserting Eqs.
(\ref{GMLS-7}) and (\ref{GMLS-8}) into Eq. (\ref{Eq-2}) yields
\begin{equation}\label{GMLS-9}
\nabla a_{j}(\bm{x})= \left( \dfrac{1}{\cos(\theta)}\dfrac{\partial a_{j}(\bm{x})}{\partial \lambda},
\dfrac{\partial a_{j}(\bm{x})}{\partial \theta} \right)^{T}:=\left(G_{\lambda},G_{\theta}\right)^{T},
\end{equation}
where
\begin{align*}
G_{\lambda}&={{{\left( {{\nabla _0}} \right)}_1}{a_j}(\bm{x})}
\left( { - \sin (\lambda )} \right) +  {{{\left( {{\nabla _0}} \right)}_2}{a_j}(\bm{x})}
\left( {\cos (\lambda)} \right),\\
\\
G_{\theta}&=
{{{\left( {{\nabla _0}} \right)}_1}{a_j}(\bm{x})}\left(
 { - \cos (\lambda )\sin (\theta )} \right) +{{{\left( {{\nabla _0}} \right)}_2}{a_j}(\bm{x})} \left
 ( { - \sin (\lambda )\sin (\theta )} \right) +  {{{\left( {{\nabla _0}} \right)}_3}{a_j}(\bm{x})}
 ( {\cos (\theta )}).
\end{align*}
\section{The MKLS formulation for advection operator}\label{Sec-4}
Our goal of this part is to introduce a new approximation namely
MKLS on the unit sphere. Previously, this technique had been given
in subdomains of $\mathbb{R}^{d}$ \cite{gu2003moving}. Here, we
first introduce the methodology of MKLS technique on the sphere, and
then \TODO{we approximate $\nabla u$ for a given function $u$
defined on the unit sphere
 in spherical coordinates using this
approach.}

We suppose that $u \in C^{m+1}(\mathbb{S}^2)$ is a function defined
on $\mathbb{S}^2$. Also, we consider a set of $N$ points on the unit
sphere. The approximation of the function $u$ by $\overline{u}$ can be given by
\begin{equation}\label{1-MK}
u({\bm{x}}) \approx \overline{u(\bm{x})} =
{\bm{Y}^T}({\bm{x}}){\bm{a}}({\bm{x}}) +
Z({\bm{x}}),\,\,\,\,\,\,\,{\bm{x}} \in \mathbb{S}^{2},
\end{equation}
where $\bm{Y}(\bm{x})$ and $\bm{a} (\bm{x})$ are the vectors of
spherical harmonics of at most $m$, and the unknown coefficient, respectively.
Also, $Z(\bm{x})$ represents the realization
of a stochastic process with mean zero, variance $\sigma^{2}$ and
non-zero covariance. The matrix form of the covariance can be written as
\begin{equation}\label{Matrix-cov}
{\mathop{\rm cov}} \left\{ {Z({{\bm{x}}_i}),Z({{\bm{x}}_j})}
\right\} = {\sigma ^2} {\bm{R}}\left[ {R({{\bm{x}}_i},{{\bm{x}}_j})}
\right],\,\,\,\,\,\,i,j = 1,2,...,N,
\end{equation}
where
 ${\bm{R}}\left[ {R({{\bm{x}}_i},{{\bm{x}}_j})} \right]$ and
${R({{\bm{x}}_i},{{\bm{x}}_j})}$  called the correlation matrix and the
correlation function between any pair of points located at
$\bm{x}_{i}$ and $\bm{x}_{j}$ on $\mathbb{S}^2$, respectively.
The following Gaussian function can be chosen as
 the correlation function
\begin{equation}\label{correlation}
R({{\bm{x}}_i},{{\bm{x}}_j}) = e^{ - c\,
\DO{dist(\bm{x}_{i},\bm{x}_{j})^2}},
\end{equation}
where $c >0$ denotes the value of the correlation parameter,
 which can be effected on the approximation solution.
 In the same manner \cite{gu2003moving},  Eq.
(\ref{1-MK}) can be written as below
\begin{equation}\label{2-MK}
\overline{u(\bm{x})} = {\bm{Y}^T}({\bm{x}}){\left(
{{{\bm{P}}^T}{{\bm{R}}^{ - 1}}{\bm{P}}} \right)^{ -
1}}{{\bm{P}}^T}{{\bm{R}}^{ - 1}}{\bf{u}} +
{\bm{r}^T}({\bm{x}}){{\bm{R}}^{ - 1}}\left( {{\bf{I}} -
{\bm{P}}{{\left( {{{\bm{P}}^T}{{\bm{R}}^{ - 1}}{\bm{P}}} \right)}^{
- 1}}{{\bm{P}}^T}{{\bm{R}}^{ - 1}}}
\right)\bm{u},\,\,\,\,\,\,\,\,\,\,{\bm{x}} \in \mathbb{S}^2,
\end{equation}
where the vector $r^{T}(\bm{x})=[R(\bm{x}_{1},\bm{x})\,\,
R(\bm{x}_{2},\bm{x}) ... R(\bm{x}_{|I(\bm{x}_{c})|},\bm{x})]$ such
that $\bm{x}_{c}$ is the evaluation point on $\mathbb{S}^2$. Due to
Eq. (\ref{2-MK}),
 the shape functions of the MKLS approximation on the unit sphere are denoted by the
 following formula
\begin{align}\label{4-MK}\nonumber
{\bm{a}^{T}}{({\bm{x}})} :=& {\bm{Y}^T}({\bm{x}}){\left(
{{{\bm{P}}^T}{{\bm{R}}^{ - 1}}{\bm{P}}} \right)^{ -
1}}{{\bm{P}}^T}{{\bm{R}}^{ - 1}} + {\bm{r}^T}({\bm{x}}){{\bm{R}}^{ -
1}}\left( {{\bf{I}} - {\bm{P}} {{\left( {{{\bm{P}}^T}{{\bm{R}}^{ -
1}}{\bm{P}}} \right)}^{ - 1}}{{\bm{P}}^T}{{\bm{R}}^{ - 1}}}
\right)\\\nonumber
\\
 =& \left[
{{a_1}({\bm{x}})\,\,\,\,\,{a_2}({\bm{x}})\,\,...\,\,\,{a_{|I({{\bm{x}}_c})|}}({\bm{x}})}
\right].
\end{align}
The approximation of a function $u$ can be given as follows
\begin{equation}\label{5-MK}
\overline{u(\bm{x})}=\displaystyle \sum_{j \in I(\bm{x})} a_{j}(\bm{x})u_{j},\,\,\,\, \bm{x} \in \mathbb{S}^2,
\end{equation}
\TODO{where $a_{j}(\bm{x}),\,\,j=1,2,...,N$ are obtained in
(\ref{4-MK}), and $u_{j}$ is the value $u$ at point $\bm{x}_{j} \in
I(\bm{x})$.} As shown in \cite{gu2003moving}, the constructed shape
functions of MKLS approximation in $\Omega \subset \mathbb{R}^d$
satisfy Kroncker's delta property. This feature also remains for the
shape functions (\ref{4-MK}) obtained on the unit sphere. Now, we
approximate the surface gradient operator $\nabla u$ using
\TODO{MKLS method}.
 The partial derivatives of $a_{j}(\bm{x})$ with respect to
$\lambda$ and $\theta$ can be obtained as
\begin{align}\label{6-MK}\nonumber
\dfrac{{\partial {a_j}(\bm{x})}}{{\partial \lambda }} &=
{{{\left( {{\nabla _0}} \right)}_1}{a_j}(\bm{x})} \dfrac{{\partial x}}{{\partial \lambda }}
+  {{{\left( {{\nabla _0}} \right)}_2}{a_j}(\bm{x})} \dfrac{{\partial y}}{{\partial \lambda }}
+  {{{\left( {{\nabla _0}} \right)}_3}{a_j}(\bm{x})} \dfrac{{\partial z}}{{\partial \lambda }}\\\nonumber
\\
  &=  {{{\left( {{\nabla _0}} \right)}_1}{a_j}(\bm{x})}
 \left( { - \sin (\lambda )\cos (\theta )} \right) + {{{\left( {{\nabla _0}} \right)}_2}{a_j}(\bm{x})}
\left( {\cos (\lambda )\cos (\theta )} \right),
\end{align}
and
\begin{align}\label{7-MK}\nonumber
\dfrac{{\partial {a_j}(\bm{x})}}{{\partial \theta }} &=
 {{{ {{\left(\nabla _0\right)}}}_1}{a_j}(\bm{x})}\dfrac{{\partial x}}{{\partial \theta }}
 + {{{\left( {{\nabla _0}} \right)}_2}{a_j}(\bm{x})} \dfrac{{\partial y}}{{\partial \theta }}
 +  {{{\left( {{\nabla _0}} \right)}_3}{a_j}(\bm{x})}\dfrac{{\partial z}}{{\partial \theta }}\\\nonumber
\\
&=  {{{\left( {{\nabla _0}} \right)}_1}{a_j}(\bm{x})}\left(
 { - \cos (\lambda )\sin (\theta )} \right) +{{{\left( {{\nabla _0}} \right)}_2}{a_j}(\bm{x})} \left
 ( { - \sin (\lambda )\sin (\theta )} \right) +  {{{\left( {{\nabla _0}} \right)}_3}{a_j}(\bm{x})}
  ({\cos (\theta )}),
\end{align}
where ${\left( {{\nabla _0}} \right)}_1$, ${\left( {{\nabla _0}}
\right)}_2$ and ${\left( {{\nabla _0}} \right)}_3$ act on the vector
functions $\bm{Y}^{T}(\bm{x})$ and ${\bm{r}^T}({\bm{x}})$ according
to Eq. (\ref{4-MK}). The partial derivatives of $\bm{Y}^{T}(\bm{x})$
with respect to $\lambda$ and $\theta$ can be obtained similar to
GMLS approximation, which is described in the previous section. On
the other hand, since ${\bm{r}^T}({\bm{x}})$ is a radial function,
its partial derivatives with respect to $\lambda$ and $\theta$ can
be computed in the same way that was given in
\cite{flyer2007transport}, and then $\nabla a_{j}(\bm{x})$ will be
computed at each point $\bm{x} \in \mathbb{S}^2$.

\section{The full-discrete scheme}\label{Sec-5}
In this section, we apply two approximations, which are given in Sections \ref{Sec-3} and \ref{Sec-4} for discretizing
the spatial variables of semi-discretized equations (\ref{BDF-3}) and (\ref{BDF-4}).
We consider $N$ points such as $X=\{\bm{x}_{1},\bm{x}_{2},...,\bm{x}_{N}\}$ on the
unit sphere in spherical coordinates, and we assume that the approximation solution
of $u^{n+1}$ at each point $\bm{x} \in
\mathbb{S}^2$ is
\begin{equation}\label{full-1}
u^{n+1}(\bm{x})\approx \displaystyle \sum_{j \in I(\bm{x})} a_{j}(\bm{x})u^{n+1}_{j},
\end{equation}
where $a_{j}(\bm{x})$ can be chosen from (\ref{GMLS-1}) or
(\ref{5-MK}), and $n=0,1,...,M-1$. The surface gradient, i.e.,
$\nabla u^{n+1}$ can be approximated by the following formula
\begin{equation}\label{full-2}
\nabla u^{n+1}(\bm{x})\approx \displaystyle
\sum_{j \in I(\bm{x})}\nabla a_{j}(\bm{x})u^{n+1}_{j},
\end{equation}
where $\nabla a_{j}(\bm{x})$
are defined due to GMLS or MKLS approximation.
Replacing Eqs. (\ref{full-1}) and (\ref{full-2}) into Eq. (\ref{BDF-4}) at each point
$\bm{x}_{i}$ for $n=0$ gives


\begin{equation}\label{full-3}
\sum_{j \in I(\bm{x}_{i})} a_{j}(\bm{x}_{i})u^{1}_{j}+\Delta t \bm{v}^{1}.\displaystyle
\sum_{j \in I(\bm{x}_{i})}\nabla a_{j}(\bm{x}_{i})u^{1}_{j}=\sum_{j \in I(\bm{x}_{i})} a_{j}(\bm{x}_{i})u^{0}_{j},
\end{equation}
where $i=1,2,...,N$. Substituting Eqs. (\ref{full-1}) and
(\ref{full-2}) into Eq. (\ref{BDF-3}) yields
\begin{equation}\label{full-4}
3\sum_{j \in I(\bm{x}_{i})} a_{j}(\bm{x}_{i})u^{n+1}_{j}+2\Delta t \bm{v}^{n+1}.\displaystyle
\sum_{j \in I(\bm{x}_{i})}\nabla a_{j}(\bm{x}_{i})u^{n+1}_{j}=4\sum_{j \in I(\bm{x}_{i})} a_{j}(\bm{x}_{i})u^{n}_{j}-
\sum_{j \in I(\bm{x}_{i})} a_{j}(\bm{x}_{i})u^{n-1}_{j},
\end{equation}
for $i=1,2,...,N$ and $n=1,2,...,M-1$.

The matrix form of Eq. (\ref{full-3}) can be written as follows
\begin{equation}\label{full-5}
\TO{A_{X}U^{1}_{X}+\Delta t
({v}^{1}_{1}.*B^1_{X})U^{1}_{X}+({v}^{1}_{2}.*B^2_{X})U^{1}_{X}=A_{X}U^{0}_{X},}
\end{equation}
\DO{where $A_{X}$, $B^1_{X}$ and $B^2_{X}$ are the global matrices
as
\[{A_X} = {\left[ {{a_j}({\bm{x}_i})} \right]
_{1 \le i \le N,1 \le j \le  N }},\,\,\,\,\,\,\,\,\,\,\,\, {B^1_X} =
{\left[ \dfrac{1}{\cos(\theta_{i})}\dfrac{ \partial
a_{j}({\bm{x}_i})}{\partial \lambda} \right]_{1 \le i \le N,1 \le j
\le N }},\]}

\[\DO{{B^2_X} =
{\left[\dfrac{ \partial a_{j}({\bm{x}_i})}{\partial \theta}
\right]_{1 \le i \le N,1 \le j \le N }}.}\]
 $U^{0}_{X}$ and
$U^{1}_{X}$ are the vectors of approximation at $t=t_{0}$ and
$t=t_{1}$, respectively. ${v}^{1}_{1}$ and $v^1_{2}$ are vectors of
velocity field at $N$ points at $t=t_{1}$, \TO{and for example in
MATLAB notation, $.*$ denotes the pointwise product between each row
of ${v}^{1}_{1}$ and each row of the matrix $B^1_{X}$.} Similarly,
Eq. (\ref{full-4}) can be represented in the following matrix form
\begin{equation}\label{full-6}
\TO{3A_{X}U^{n+1}_{X}+2\Delta
t(({v}^{n+1}_{1}.*B^1_{X})U^{n+1}_{X}+({v}^{n+1}_{2}.*B^2_{X})U^{n+1}_{X})=4A_{X}U^{n}_{X}-A_{X}U^{n-1}_{X},}
\end{equation}
where $U^{n-1}_{X}$, $U^{n}_{X}$ and $U^{n+1}_{X}$ are the vectors
of approximation at $t=t_{n-1}$, $t=t_{n}$
and $t=t_{n+1}$, respectively.

\DO{In order to solve the linear system of algebraic equations
    obtained here, i.e., (\ref{full-5}) and (\ref{full-6}), an iterative
    algorithm namely the biconjugate gradient stabilized (BiCGSTAB) method with zero-fill incomplete lower-upper (ILU)
     preconditioner is employed. In the
    literature \cite{lehto2017radial}, it has been shown that the method is efficiently
    solvable for the linear system of algebraic equations generated by
    the third-order semi-implicit backward differential formula and
    combined with a meshless technique for the solution of
    reaction-diffusion equation on the surfaces \cite{lehto2017radial}.
    It also should be noted that this algorithm can be used only for the
    linear system of algebraic equations with a large sparse coefficient
    matrix \cite{lehto2017radial}. According to the approximations
    presented here, the final coefficient matrices in Eqs.
    (\ref{full-5}) and (\ref{full-6}) are sparse, and thus the BiCGSTAB
    algorithm could be applied without any difficulty as we will observe in
    the next section.}

\begin{figure}[ht]
    \centering
    \includegraphics[width=16cm,height=12cm]{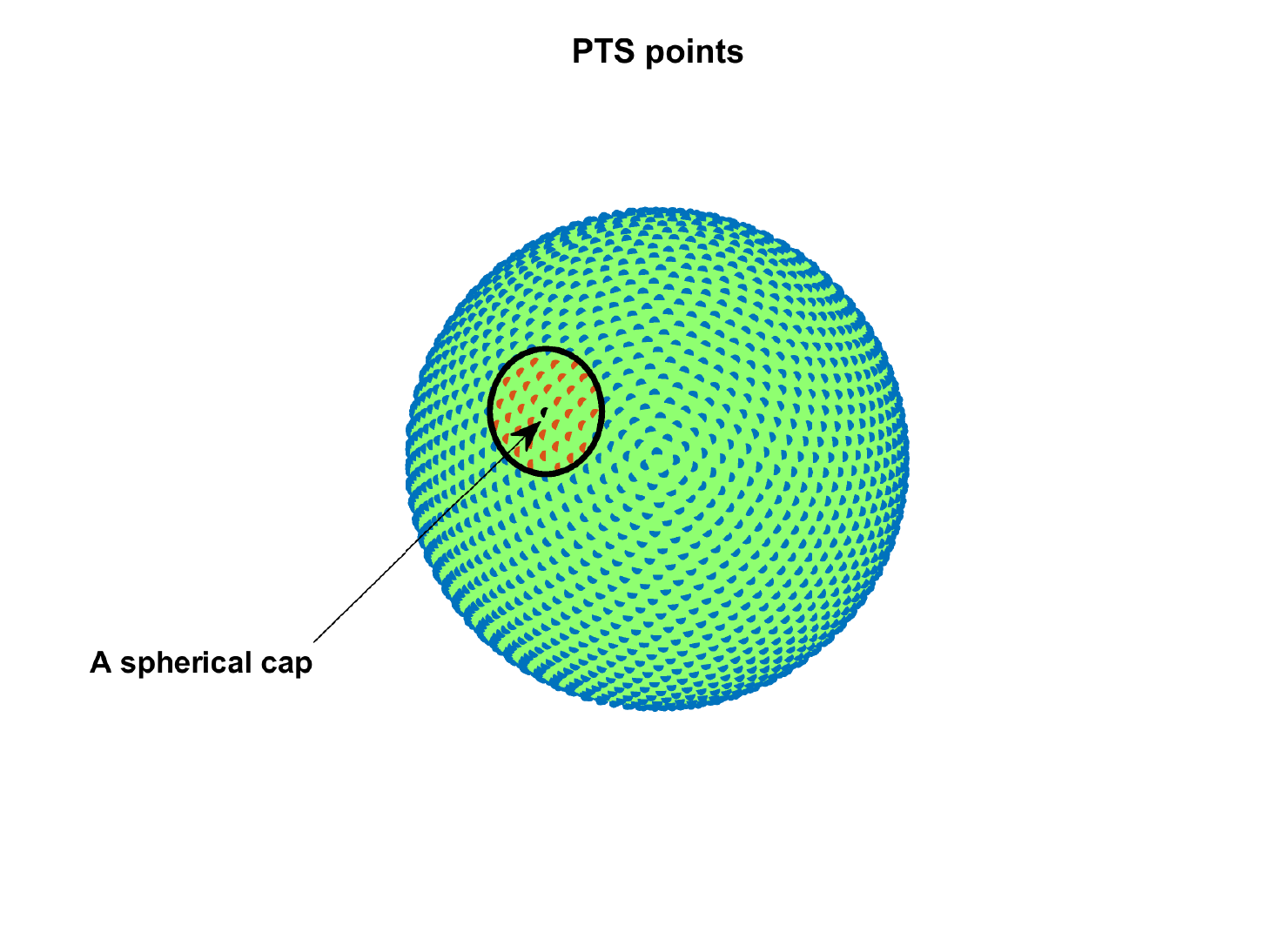}
    \vspace{-2cm}
    \caption{Set of distribution points on the sphere for PTS
        points with an example of spherical cap of radius
        $\delta$.}
    \label{fig-0}
\end{figure}
\section{Numerical results}\label{Sec-6}
In this section, in order to investigate the ability of the proposed
method we use the three standard tests, which have been proposed in the
literature
\cite{fornberg2011stabilization,nair1999cascade,nair2010class,shankar2018mesh}.
The first case is the known "solid-body rotation of a cosine bell"
\cite{shankar2018mesh}, the second one is called "vortex roll-up"
\cite{fornberg2011stabilization,nair1999cascade}, and the latest is
"deformational flow" \cite{nair2010class}. The numerical results
reported here via quasi-uniformly distributed point sets $X$, which
are known as the \TODO{minimum energy (ME)
 \cite{flyer2007transport,womersley2003interpolation} and
phyllotaxis spiral (PTS) \cite{shankar2018mesh}}, and their fill
distance $h$ is proportional to $N^{-1/2}$, where $N$ is the number
of points distributed on the unit spheres. For the readers convenience, the
procedure of the presented numerical methods is discussed in Algorithm 1 and Algorithm 2.
\begin{algorithm} \caption{\TO{Computational algorithm of GMLS (or
MKLS) approximation}}\label{Algorithm1}
\begin{algorithmic}\vspace{0.15cm}
\State \textbf{Input:} data points on $\mathbb{S}^2$,
$X=[\bm{x}_{1},\bm{x}_{2},...,\bm{x}_{N}]^{T} \in \mathbb{R}^{N
\times 3}$, $\delta>0$ as a radius of local spherical cap; 
\vspace{0.15cm}\State
\textbf{Output:} Matrices $A_{X}$, $B^{1}_{X}$, $B^{2}_{X}$; \vspace{0.15cm}\For
{$i=1,2,...,N$}
\vspace{0.15cm}\State Construct $I(\bm{x}_{i})$ due to $X$;
\vspace{0.15cm}\State Compute local matrix $P^{T}WP$ due to points in
 $I(\bm{x}_{i})$;
\vspace{0.15cm}\State Compute the vector $Y(\bm{x})$, $\nabla_{0}Y(\bm{x})$ defined in Eqs. (\ref{GMLS-3}) and
 (\ref{GMLS-6}), respectively at each point $\bm{x}_{i}$;
\vspace{0.15cm}\State Compute the vectors $\bm{a}^{\star}(\bm{x})$ and
 $\bm{a}^{\star}_{\nabla_{0}}(\bm{x})$ at each point $\bm{x}_{i}$
 due to Eqs. (\ref{GMLS-3}) and (\ref{GMLS-6})

 (or similarly compute Eq. (\ref{4-MK}) and components of
 $\bm{a}^{\star}_{\nabla_{0}}(\bm{x})$)
  for MKLS approximation );
\vspace{0.15cm}\State $A_{X}(i,:) \gets \bm{a}^{\star}(\bm{x}_{i})$;
\vspace{0.15cm}\State $B^{1}_{X}(i,:) \gets \bm{a}^{\star}_{(\nabla_{0})_{1}}(\bm{x}_{i})$;
\vspace{0.15cm}\State $B^{2}_{X}(i,:) \gets \bm{a}^{\star}_{(\nabla_{0})_{2}}(\bm{x}_{i})$;
\vspace{0.15cm} \EndFor
\vspace{0.15cm}
\end{algorithmic}
\end{algorithm}

\begin{algorithm}
\caption{\TO{Computational algorithm for solving transport equation
on $\mathbb{S}^2$}}\label{Algorithm2}
\begin{algorithmic}
	\vspace{0.15cm}
\State \textbf{Input:} data points on $\mathbb{S}^2$,
$X=[\bm{x}_{1},\bm{x}_{2},...,\bm{x}_{N}]^{T} \in \mathbb{R}^{N
\times 3}$, $\delta>0$ as a radius of local spherical cap; $T$ as
final time, time step $\Delta t$;
 \vspace{0.15cm}\State \textbf{Output:} Approximation solution $U_{X}$ at final time;
 \vspace{0.15cm}\State Calling $A_{X},B^{1}_{X},B^{2}_{X}$ from Algorithm
 \ref{Algorithm1};
\vspace{0.15cm}\State Set the initial condition
$U^{0}_{X}=\{u(\bm{x}_{i},0)\}_{i=1}^{N},t=0,m=0$;
\vspace{0.15cm} \State Set $m=1$, $t=m \Delta t$ and compute
the velocity vector i.e., $\bm{v}$ at this time;
\vspace{0.15cm} \State Use zero-fill incomplete lower-upper (ILU)
     preconditioner for Eq. (\ref{full-5});
   \vspace{0.15cm}  \State Find $U^1_{X}$ by solving the linear system
     (\ref{full-5});
     \While {$t \leq T$}
   \vspace{0.15cm}   \State Set $m=m+1$, $t=m \Delta t$ and compute
      the velocity vector i.e., $\bm{v}$ at this time;
   \vspace{0.15cm}    \State Use zero-fill incomplete lower-upper (ILU)
     preconditioner for Eq. (\ref{full-6});
      \vspace{0.15cm}   \State Find $U^m_{X}$ by solving the linear system
     (\ref{full-6});
     \EndWhile \vspace{0.15cm}
\end{algorithmic}
\end{algorithm}

Figure \ref{fig-0} illustrates $2500$ PTS points with a spherical
cap. The $\ell_{2}$ norm is computed to show the accuracy of the
proposed two meshless methods that is defined as follows
\begin{equation}\label{Integration-1}
\left(\displaystyle
\int_{\mathbb{S}^2}[f(\bm{x})]^{2}d\bm{x}\right)^{\frac{1}{2}}
\approx \left(\displaystyle \dfrac{4 \pi}{N}\sum_{j=1}^{N}
[f(\bm{\eta}_{j})]^{2}\right)^{\frac{1}{2}}:=\|f\|_{\ell_{2}},
\end{equation}
 where $\{\bm{\eta}_{1},\bm{\eta}_{2},...,\bm{\eta}_{N}\}$ is
a set of $N$ spherical $t$-design points on the unit sphere
\cite{atkinson2012spherical,womersley2003interpolation}.
 All simulations presented here are run on a $2.2$ GHz
Intel Core i7-2670QM CPU and $8$ GB of RAM and all self-developed
codes are written \TODO{in MATLAB (version $2017$a) in standard
double precision.}

\begin{figure}[t!]
    \centering
    \includegraphics[width=6.25cm,height=5.25cm]{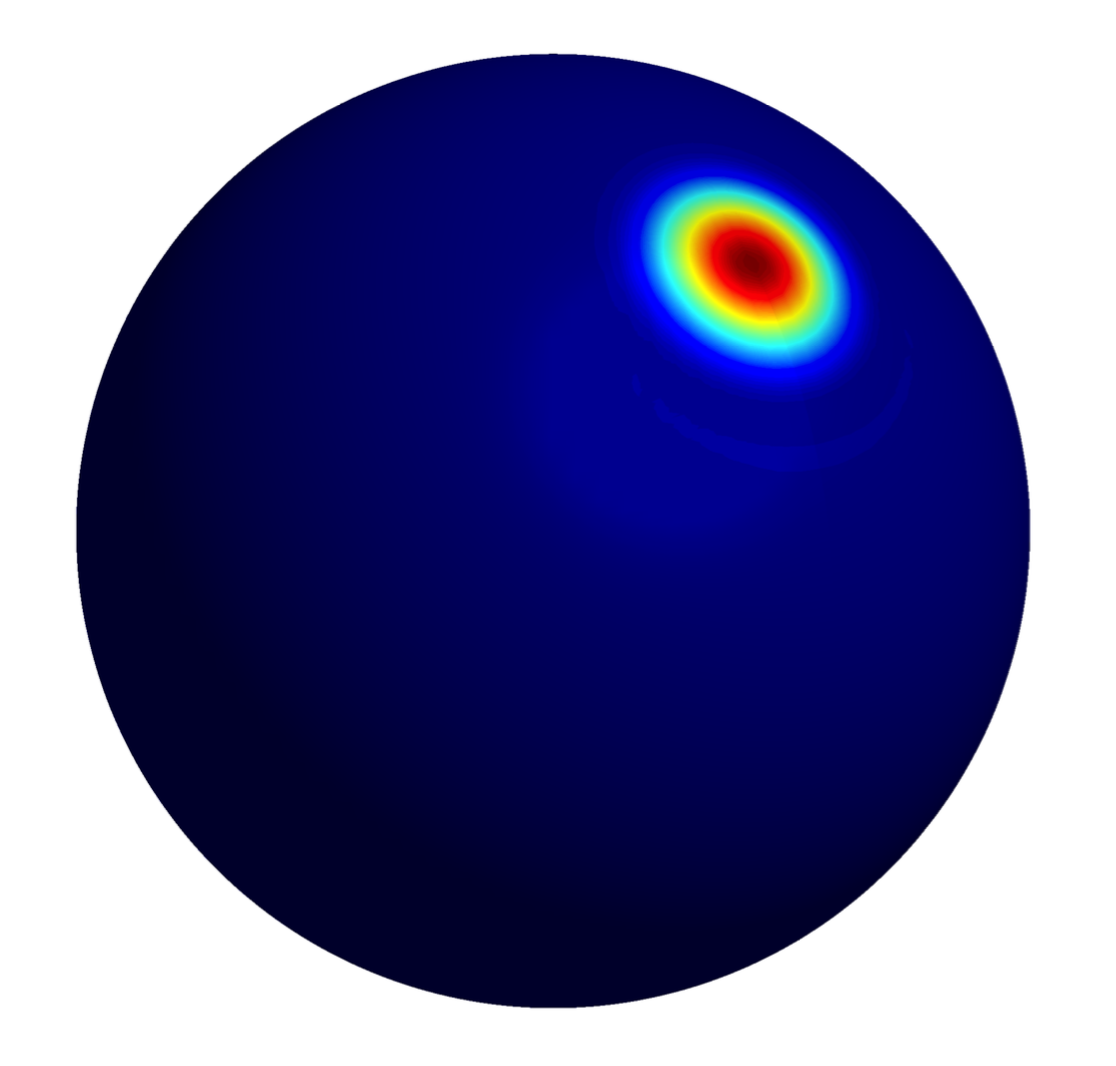}\hspace{1cm}
    \includegraphics[width=6.25cm,height=5.25cm]{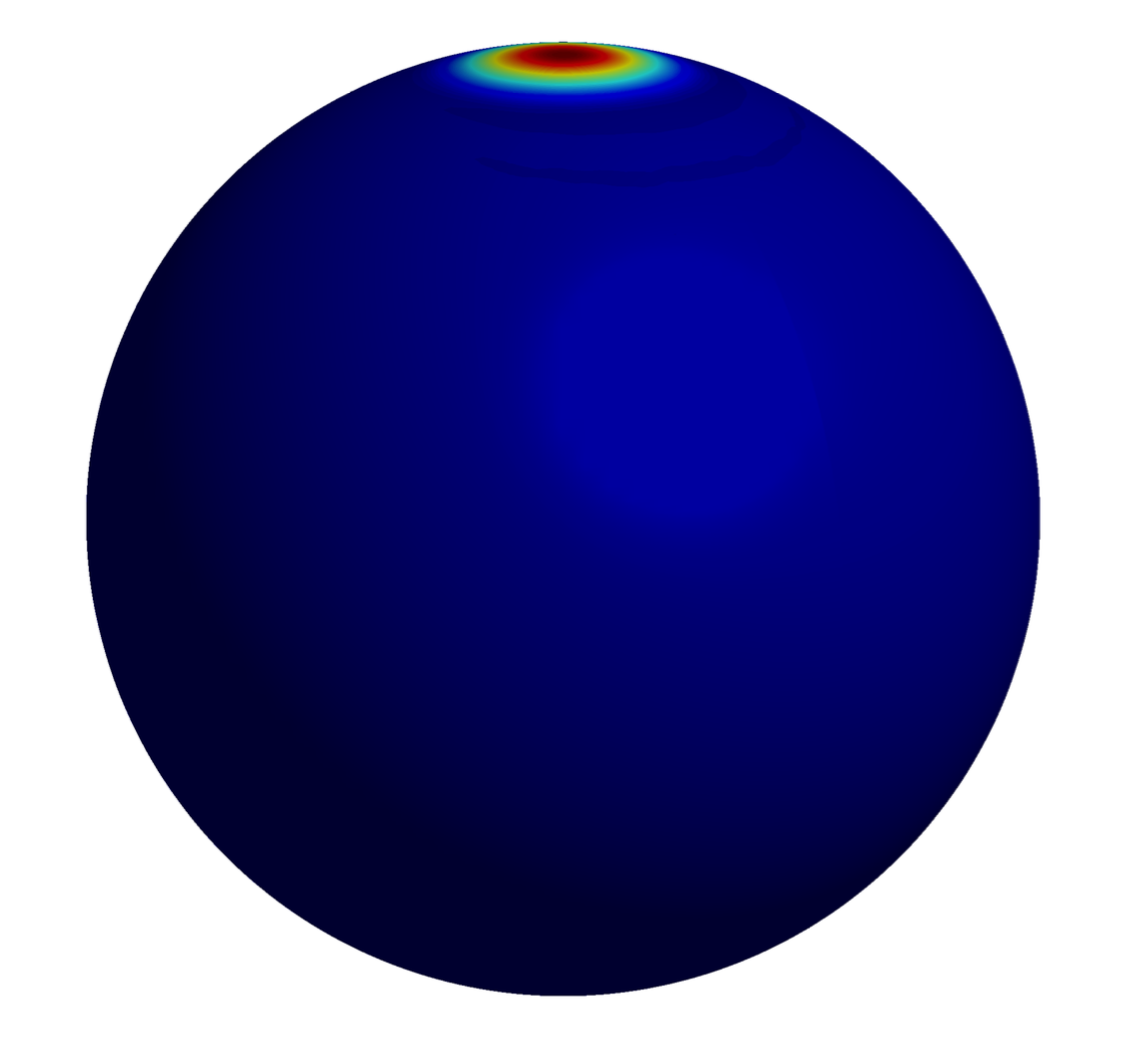}\vspace{1cm}
    \includegraphics[width=6.25cm,height=5.25cm]{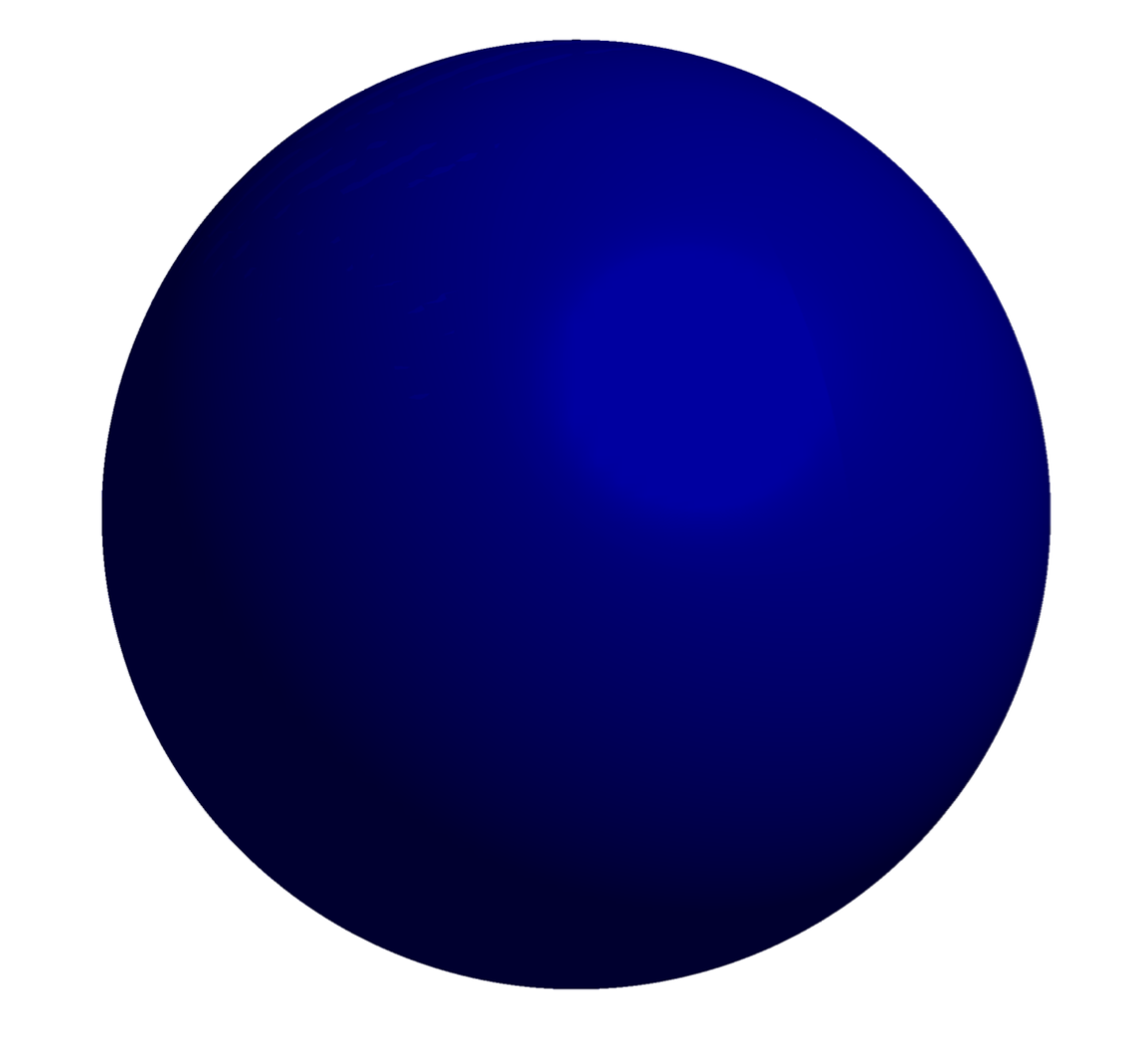} \hspace{1cm}
    \includegraphics[width=6.25cm,height=5.5cm]{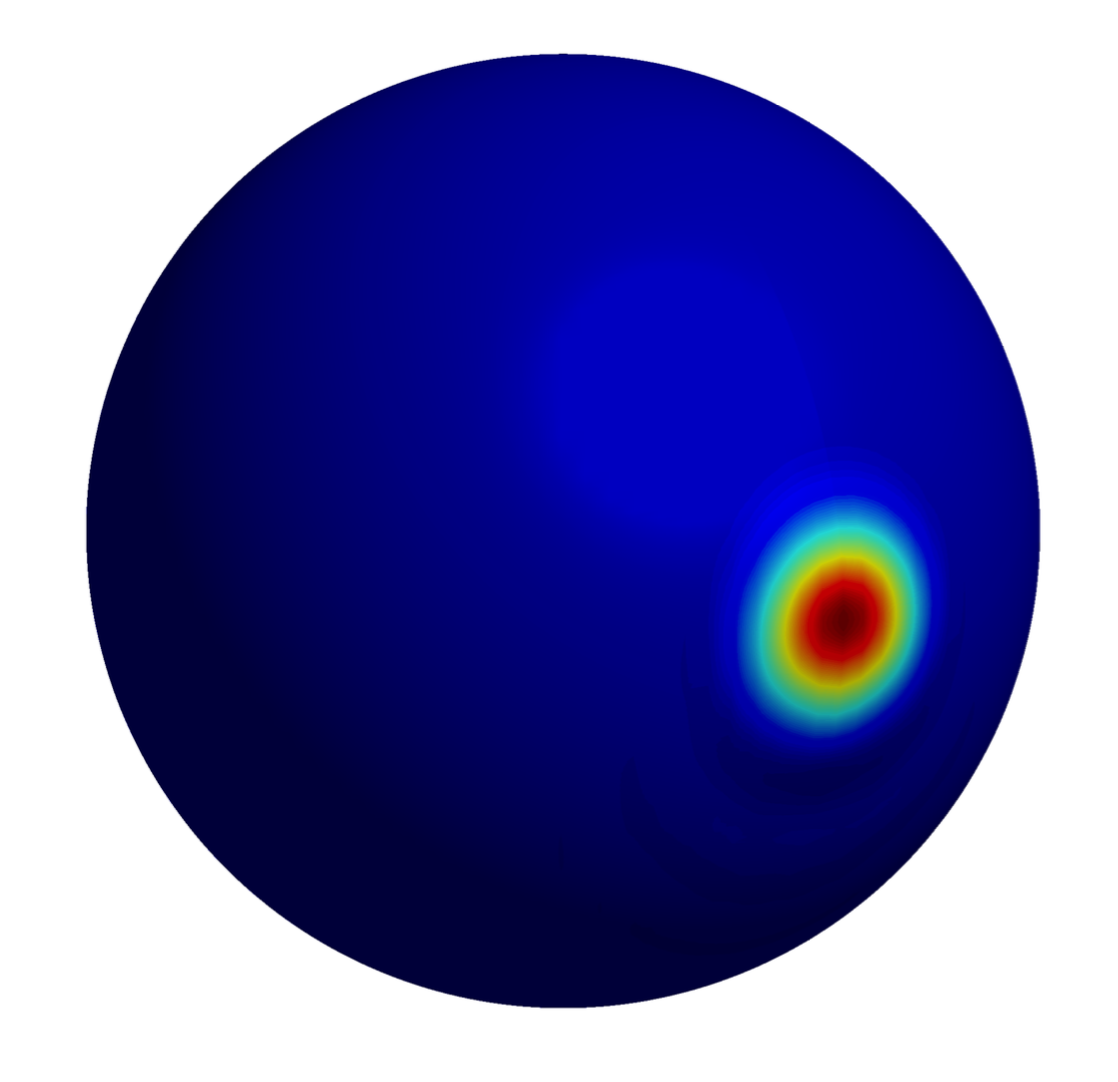}
    \caption{The 3D simulation of one full revolution of the bell over the sphere via GMLS approximation
        for solid-body rotation test (given in Subsection \ref{61}) at $t=T/8$ (top left),  $t=T/4$ (top right),  $t=T/2$ (bottom left), and  $t=T$ (bottom right).}
    \label{fig-1}
\end{figure}

\begin{figure}[t!]
    \centering
    \includegraphics[width=5.8cm,height=5.25cm]{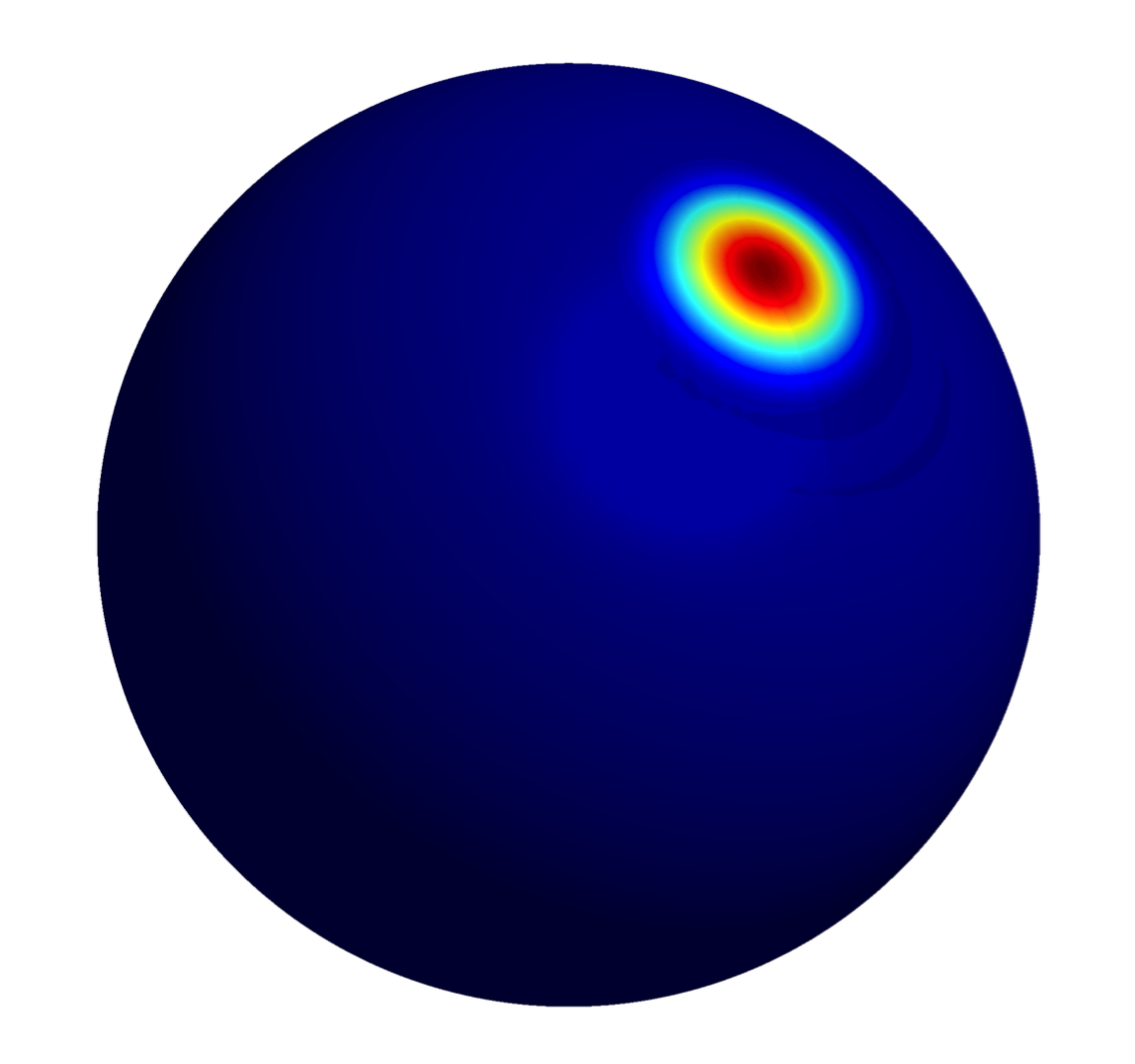}\hspace{1cm}
    \includegraphics[width=6.05cm,height=5.25cm]{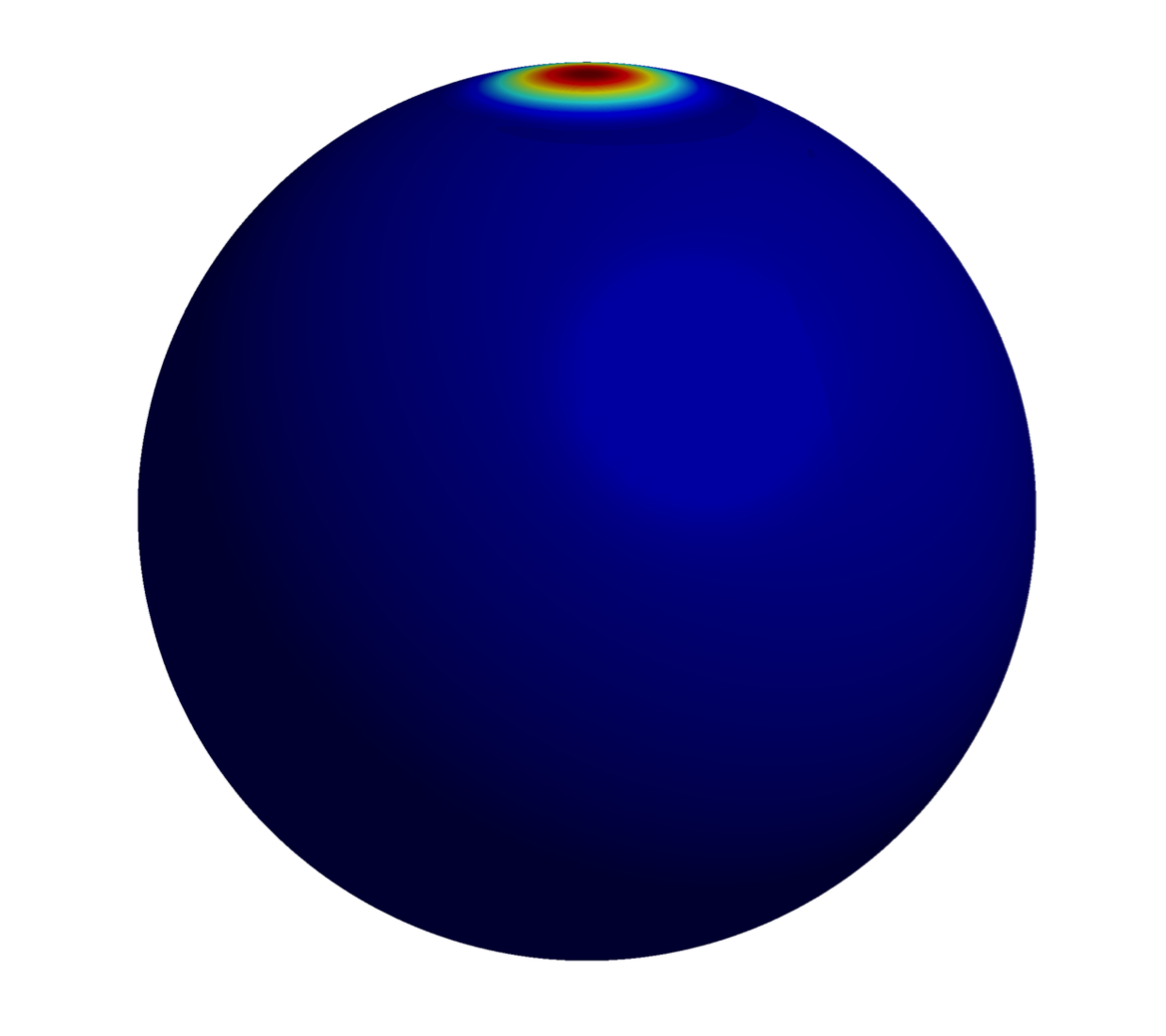}\vspace{1cm}
    \includegraphics[width=6.05cm,height=5.25cm]{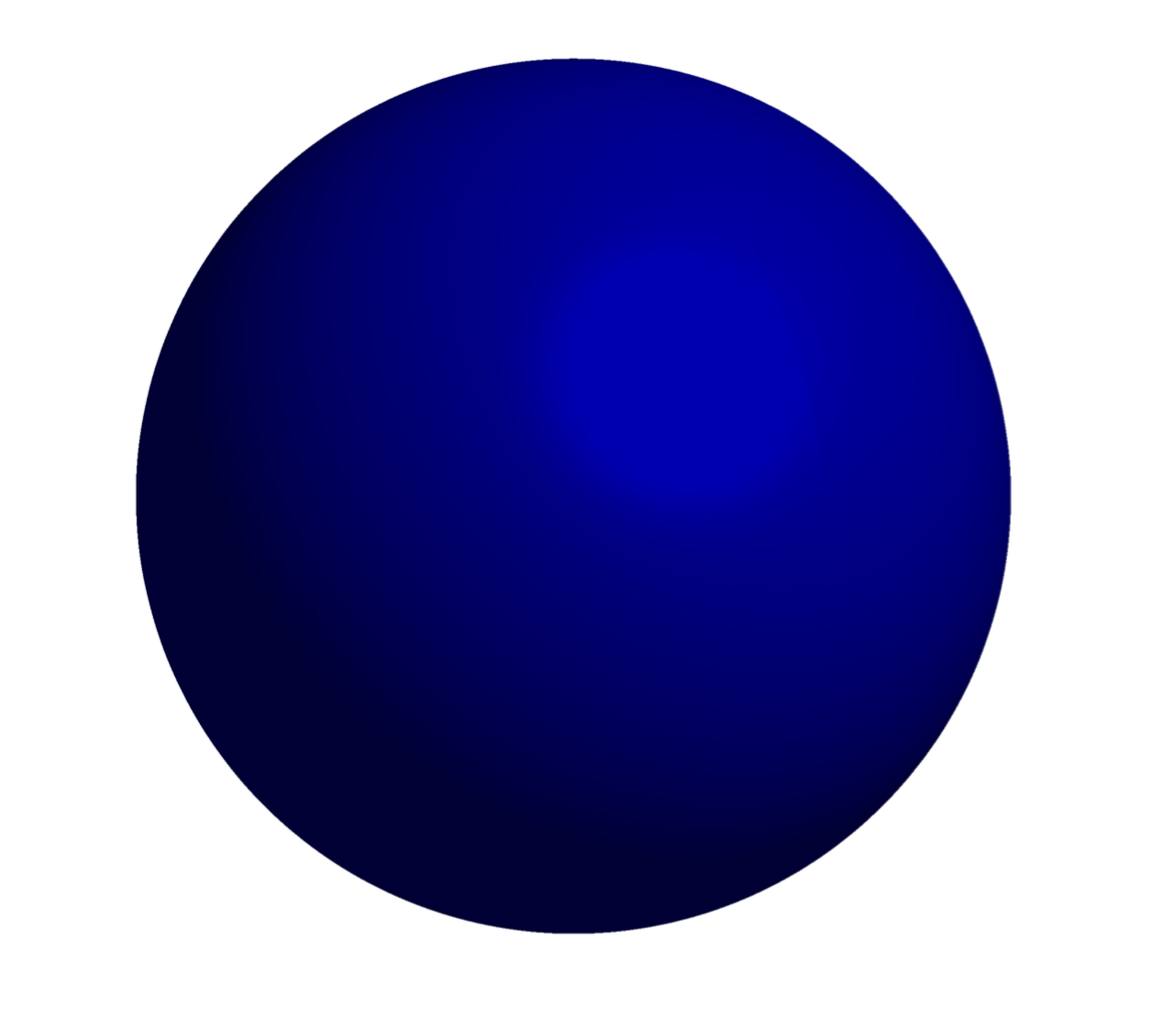}\hspace{1cm}
    \includegraphics[width=6.05cm,height=5.25cm]{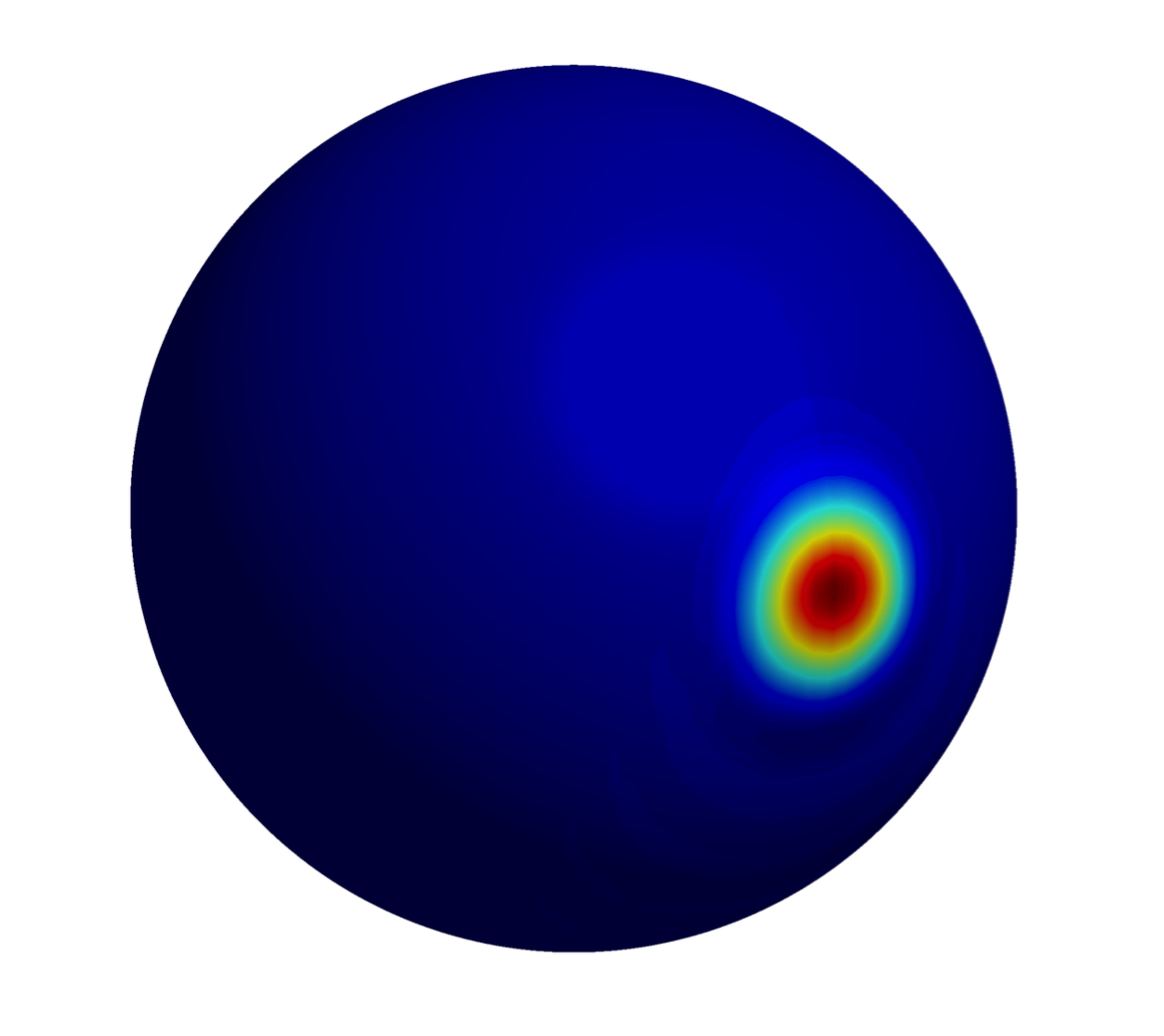}
    \caption{The 3D simulation of one full revolution of the bell over the sphere via MKLS approximation
        for solid-body rotation test (given in Subsection \ref{61}) at $t=T/8$ (top left),
        $t=T/4$ (top right),  $t=T/2$ (bottom left), and  $t=T$ (bottom right).}
    \label{fig-2}
\end{figure}

\begin{figure}[t!]
    \centering
    \includegraphics[width=8.25cm,height=5.5cm]{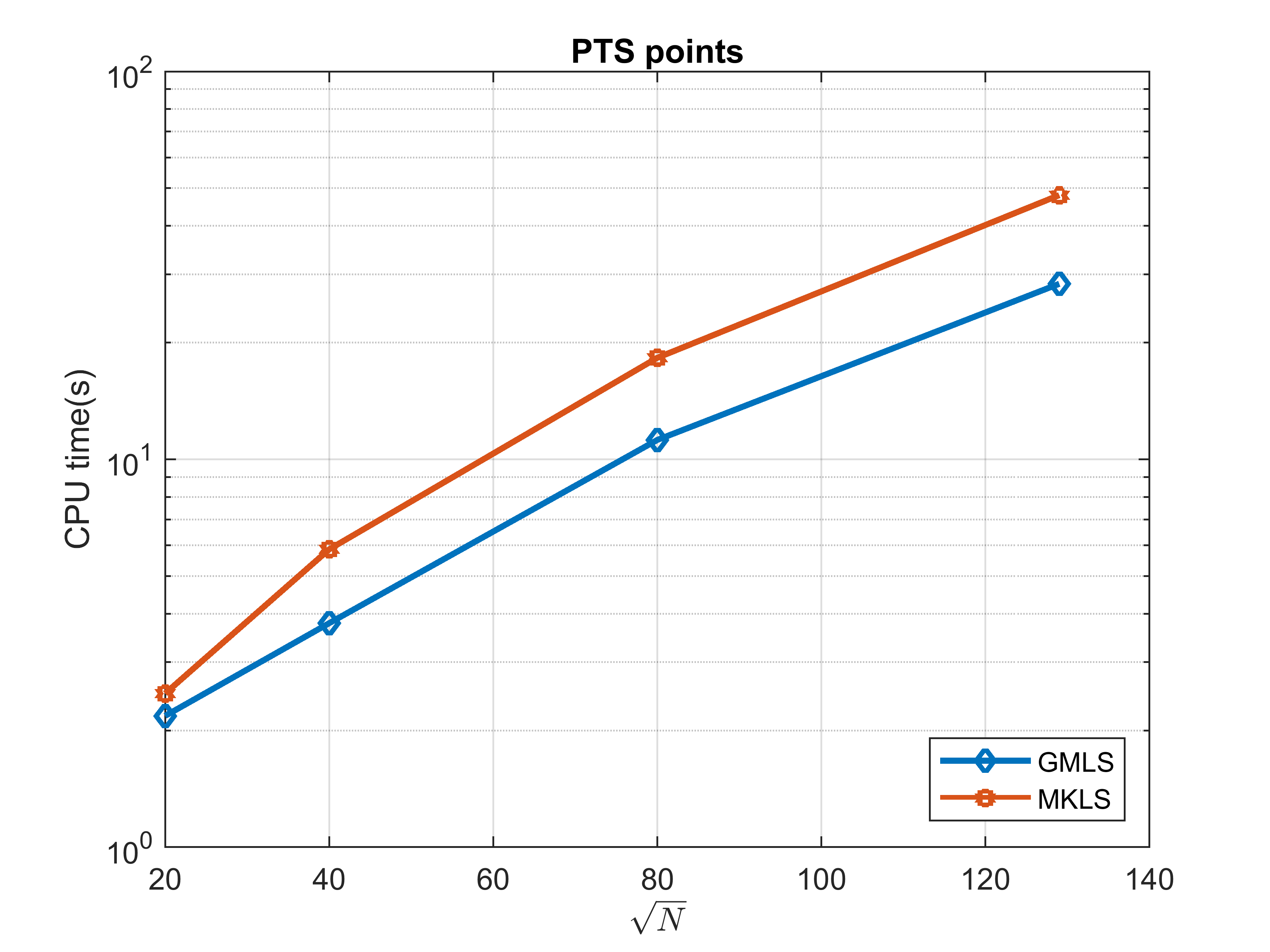}\hspace{1cm}
    \includegraphics[width=8.25cm,height=5.5cm]{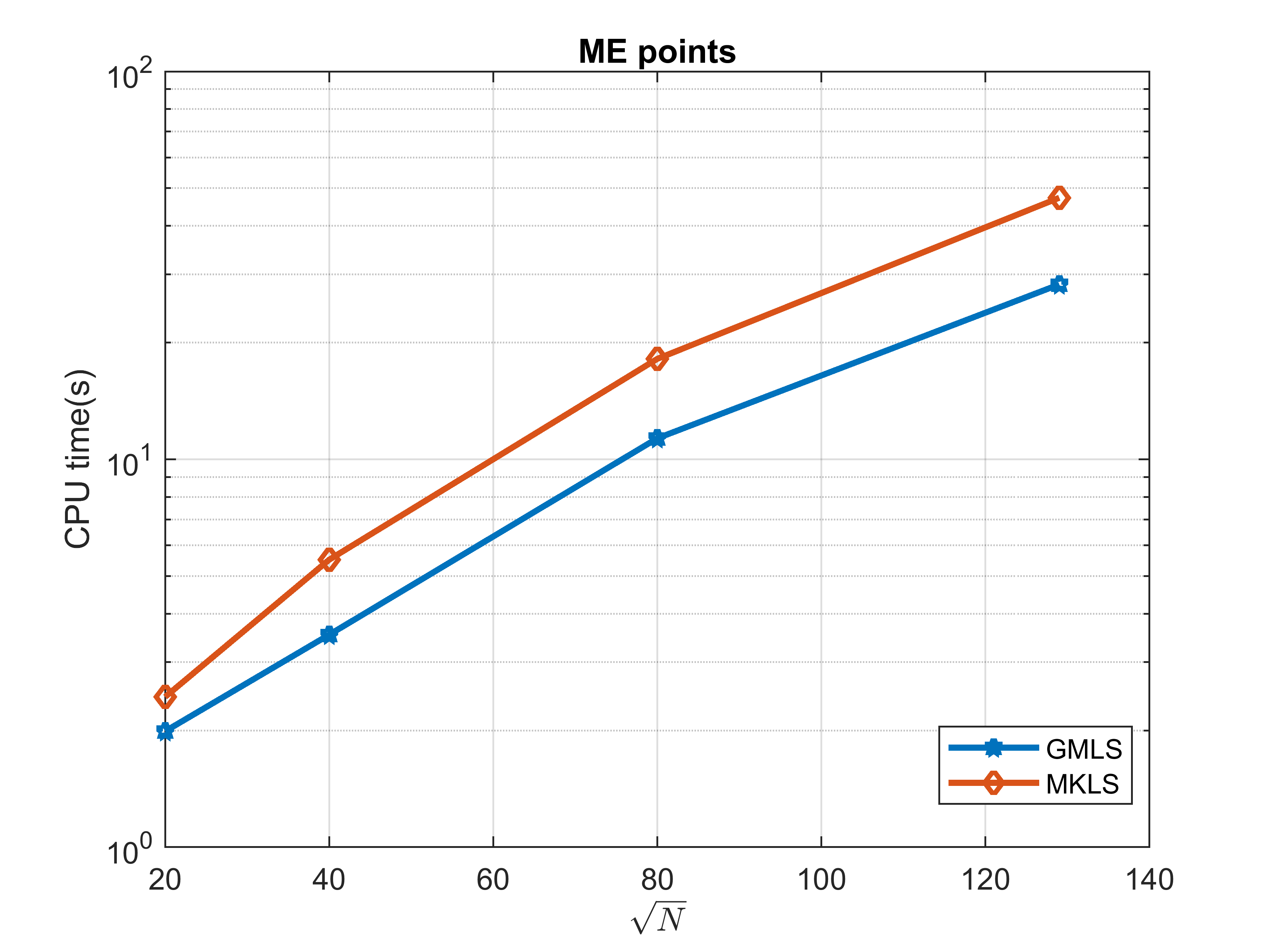} 
    \caption{\TO{The used CPU time for constructing all required matrices in
    Algorithm \ref{Algorithm1} in GMLS and MKLS approximations using
    different values of $N$, PTS points (left) and ME points (right).}}
    \label{fig-2-1}
\end{figure}

\subsection{Solid-body rotation of a cosine bell test}\label{61}
As the first standard test, we considered the transport equation (\ref{Eq-1})
on the unit sphere
with the following
vector velocity field \cite{williamson1992standard,shankar2018mesh}
$$v_{1}(\lambda,\theta)=\sin(\theta)\sin(\lambda)\sin(\alpha)-\cos(\theta)\cos(\alpha),\,\,\,\,\,\,
v_{2}(\lambda,\theta)=\cos(\lambda)\sin(\alpha),$$
where $-\pi \leq \lambda \leq \pi$ and $-\pi/2 \leq \theta \leq \pi/2$. Here, we have chosen $\alpha=\pi/2$, which
 shows the advection of the initial condition over the north and
south poles directly \cite{shankar2018mesh}.
The initial condition for this test is considered as follows \cite{shankar2018mesh}
\begin{equation} \label{Initial-1}
u(\lambda ,\theta,t=0 ) = \left\{ \begin{array}{l}
\dfrac{1}{2}\left( {1 + \cos \left( {\dfrac{{\pi r}}{{{R_b}}}} \right)} \right),\,\,\,\,\,\,\,\,\,\,r < {R_b},\\
0,\,\,\,\,\,\,\,\,\,\,\,\,\,\,\,\,\,\,\,\,\,\,\,\,\,\,\,\,\,\,\,\,\,\,\,\,\,\,\,\,\,\,\,\,\,
\,\,\,\,\,\,\,\,\,\,\,\,\,r \ge {R_b}.
\end{array} \right.
\end{equation}
Here $r=\arccos(\cos(\theta)\cos(\lambda))$ and $R_{b}=\frac{1}{2}$.
This example illustrates one full revolution of the bell over the
sphere at $T=2 \pi$. To simulate this process via the meshless
methods presented here, we have fixed $N=19600$ PTS points,
 $\delta=12h$ with $h=N^{-1/2}$
and $\Delta t=T/1000$. Also, the constant parameter in MKLS method
is considered experimentally $c=20/h$ \cite{gu2003moving}. Figures
\ref{fig-1} and \ref{fig-2} show the numerical solutions
 of $u$ at different time levels $t=T/8,~T/4,~T/2$ and
$t=T$ using GMLS and MKLS approximations. \TODO{The results obtained
via two methods are in a good agreement with those reported in the
literature \cite{williamson1992standard,shankar2018mesh}.} \TO{In
Figure \ref{fig-2-1}, the used CPU time for constructing all
required matrices in Algorithm \ref{Algorithm1} for both
approximations are given using different values of $N$.} \DO{Table
\ref{Table1-1} shows the used CPU time in both techniques for
different values $N$ during the above simulations.} In Tables
\ref{Table-1} and \ref{Table-2}, $\ell_{2}$ errors are computed for
different values $N$ via both techniques, respectively. \TODO{ As
can be observed in results, GMLS and MKLS approximation have almost
the same accuracy in solving this example.}

\begin{table}
\begin{center}
\begin{tabular}{lllllllllllllllllll}
  \hline
  $\textbf{Method}$&&&&$N$ &&&&$\textbf{CPU time}\,(s)$  \\ 
  \hline
\vspace{0.1cm} $\textbf{GMLS}$&&&&$1600$&&&& $3.21$  \\ \vspace{0.1cm}
 &&&&$6400$&&&& $10.66$   \\ \vspace{0.1cm}
 &&&&$19600$&&&& $28.73$  \\
 \hline 
\vspace{0.1cm}$\textbf{MKLS}$&&&&$1600$&&&& $3.26$  \\
 &&&&$6400$&&&& $10.96$   \\\vspace{0.1cm}
 &&&&$19600$&&&& $26.50$   \\
  \hline
\end{tabular}
\caption{\DO{The used CPU time  with different values $N$ \\for the BiCGSTAB method
for the first test.}}\label{Table1-1}
\end{center}
\end{table}

\begin{table}
\begin{center}
\begin{tabular}{llllllllllllll}
  \hline
  &&\multicolumn{2}{l}{PTS}&&\multicolumn{2}{l}{ME} \\
  \cline{3-4} \cline{6-7}
  $N$&&$\ell_{2}$ &&&& $\ell_{2}$ \\ 
  \hline
  \vspace{0.1cm}$400$  && $2.59\e-1$     &&&& $2.53\e-1$      \\  \vspace{0.1cm}
 $1600$  && $1.72\e-1$    &&&& $1.71\e-1$    \\  \vspace{0.1cm}
 $6400$ && $4.66\e-2$   &&&& $4.64\e-2$    \\  \vspace{0.1cm}
 $16641$ && $2.05\e-2$  &&&& $2.06\e-2$  \\
  \hline
\end{tabular}
\caption{The $\ell_{2}-$error for different values $N$ \\
 in GMLS approximation for the first test problem.}\label{Table-1}
\end{center}
\end{table}

\begin{table}
\begin{center}
\begin{tabular}{llllllllllllll}
  \hline
  &&\multicolumn{2}{l}{PTS}&&\multicolumn{2}{l}{ME} \\
  \cline{3-4} \cline{6-7}
  $N$&&$\ell_{2}$ &&&& $\ell_{2}$ \\
  \hline
 \vspace{0.1cm} $400$  && $2.68\e-1$     &&&& $2.55\e-1$      \\  \vspace{0.1cm}
 $1600$  && $2.47\e-1$    &&&& $2.43\e-1$    \\  \vspace{0.1cm}
 $6400$ && $9.53\e-2$   &&&& $8.73\e-2$    \\  \vspace{0.1cm}
 $16641$ && $2.36\e-2$  &&&& $1.97\e-2$  \\ 
  \hline
\end{tabular}
\caption{The $\ell_{2}-$error for different values $N$ \\
in MKLS approximation for the first test.}\label{Table-2}
\end{center}
\end{table}

\begin{figure}[t!]
    \centering
    \includegraphics[width=6.5cm,height=5.25cm]{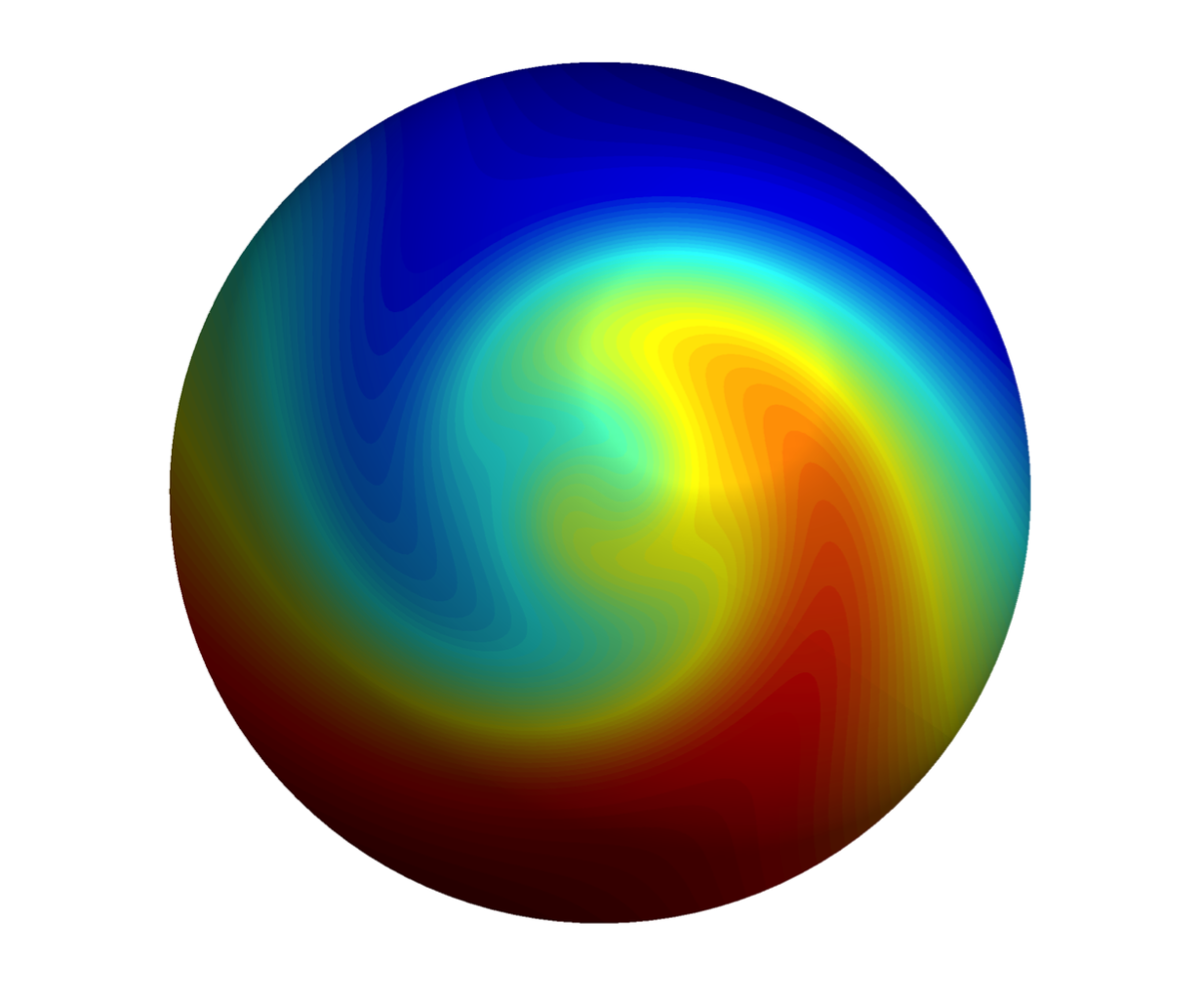}\hspace{1.5cm}
    \includegraphics[width=6.5cm,height=5.25cm]{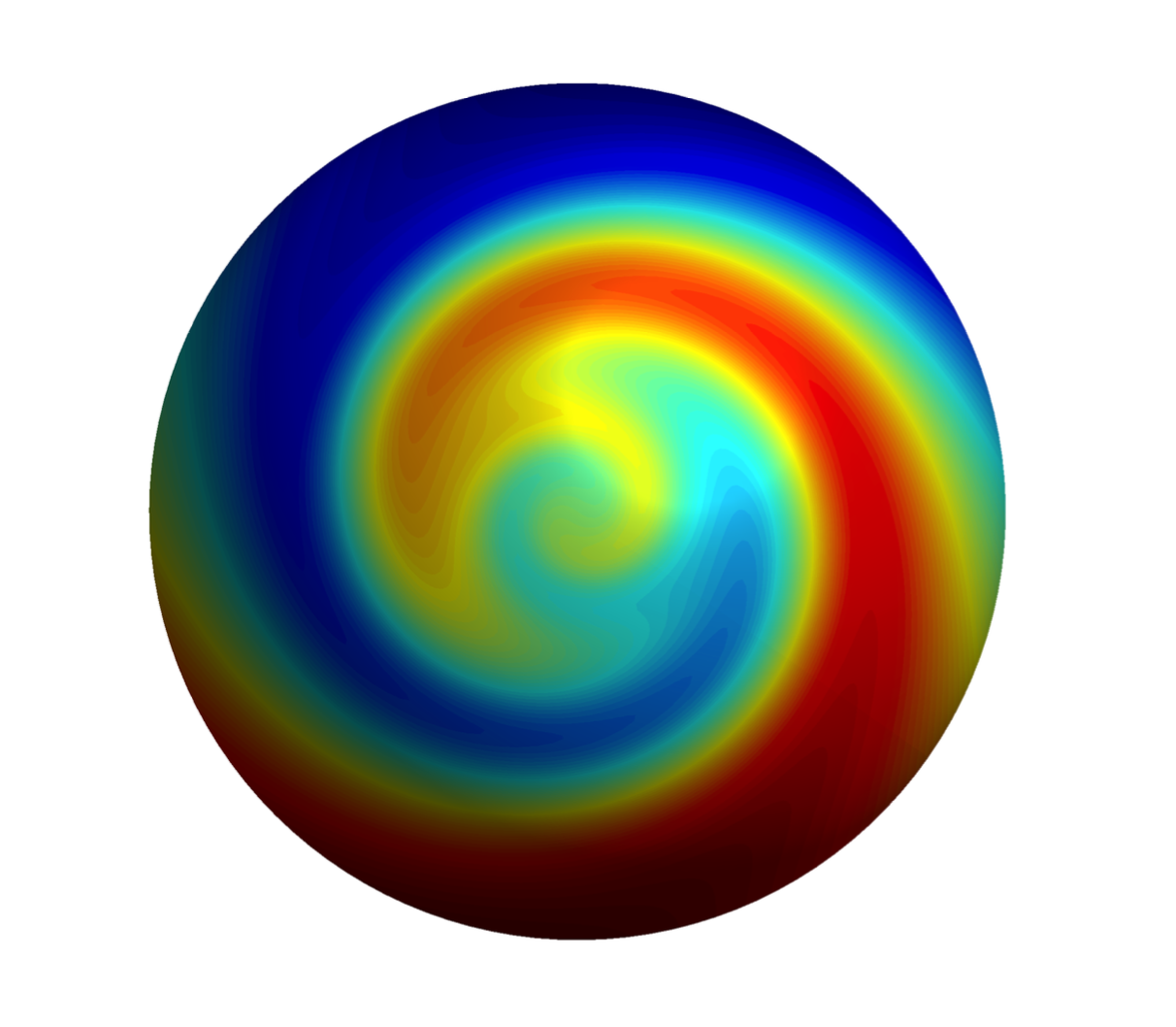}
    \includegraphics[width=6.5cm,height=5.25cm]{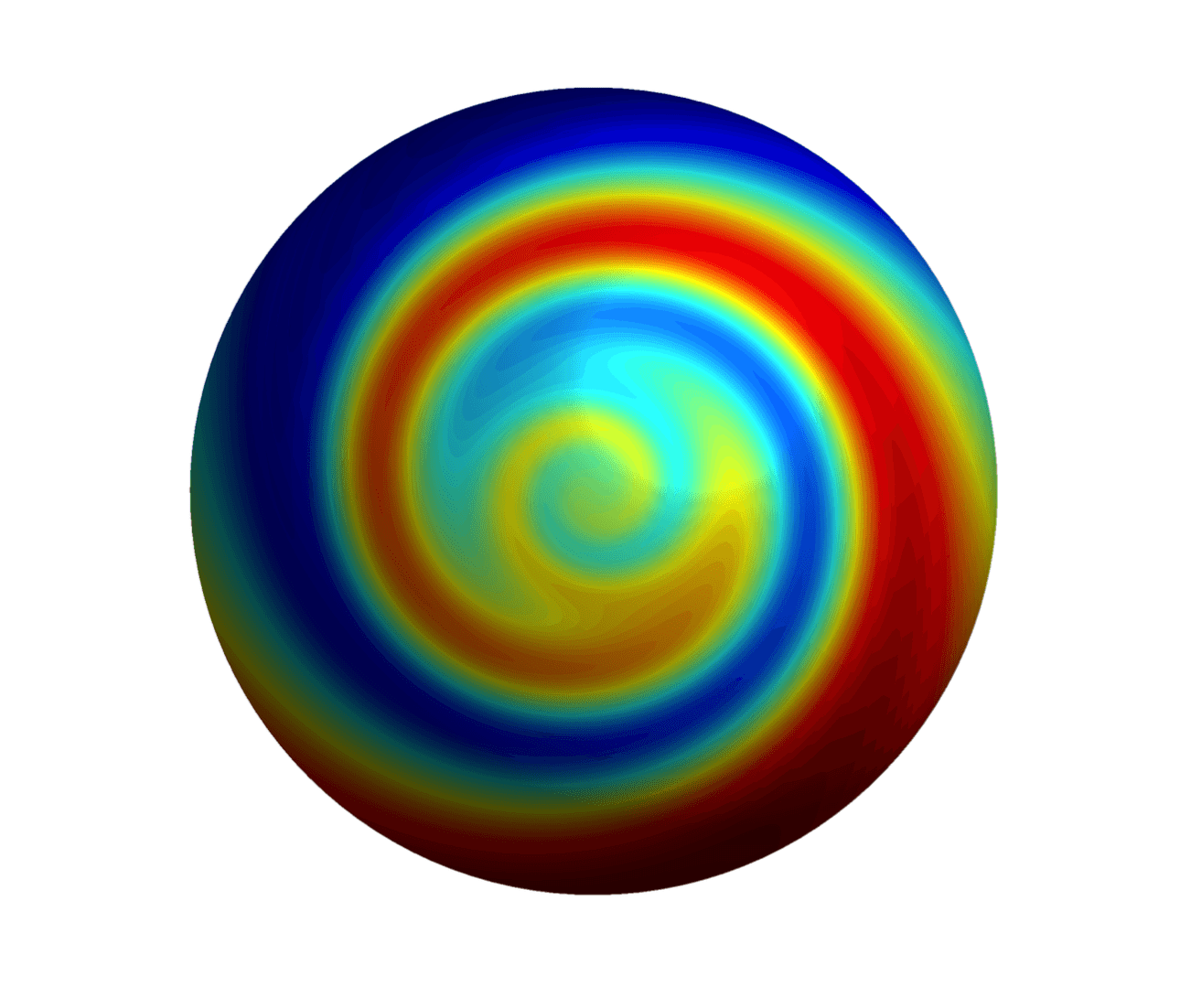}
    \caption{The 3D simulation of the vortex roll-up via GMLS approximation
        for the second test (given in Subsection \ref{62}) at $t=3$ (top left), $t=6$ (top right), and $t=9$ (bottom).}
    \label{fig-3}
\end{figure}
\subsection{Vortex roll-up test}\label{62}
\begin{figure}[t!]
    \centering
    \includegraphics[width=6.5cm,height=5.25cm]{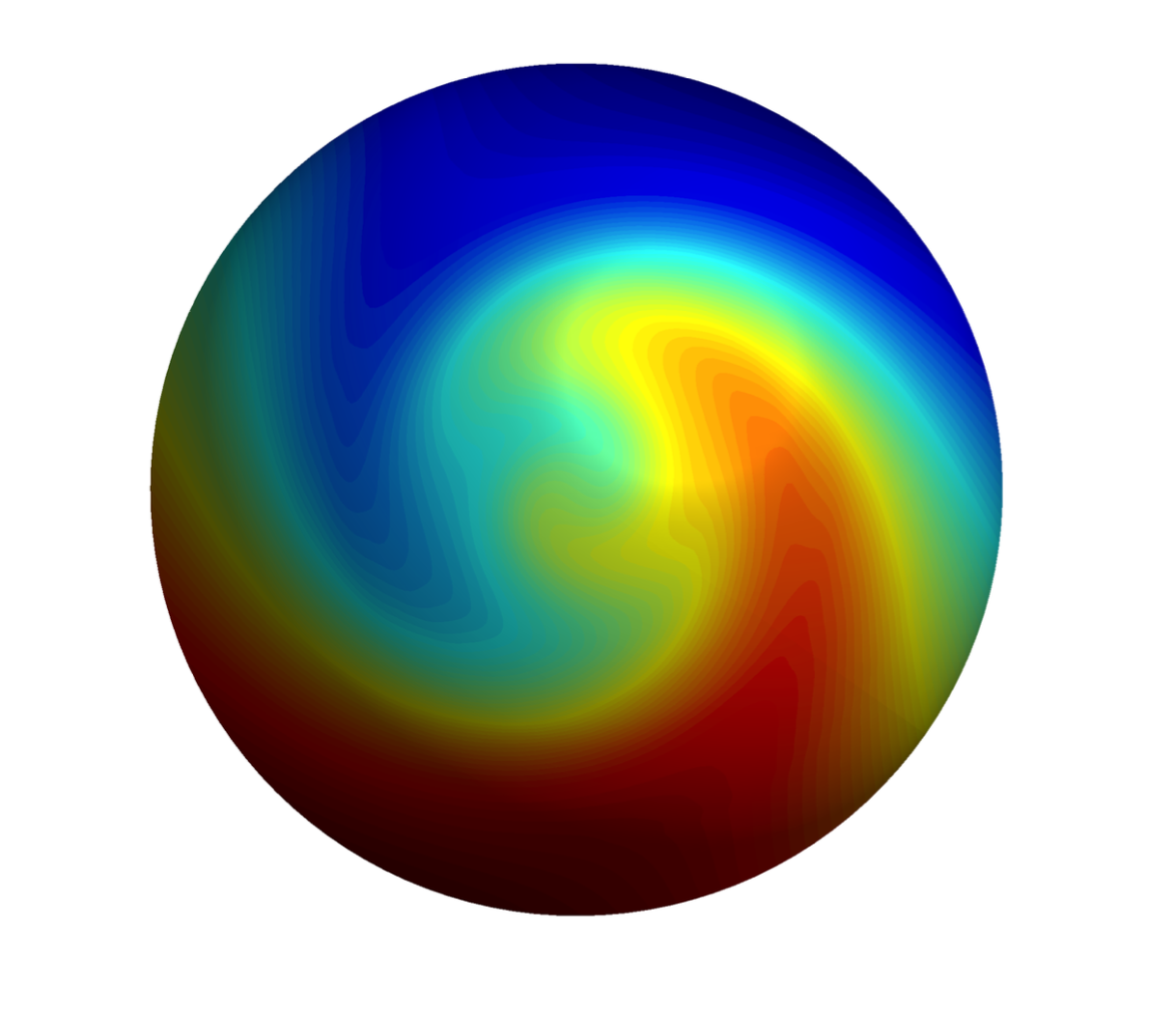}\hspace{1.5cm}
    \includegraphics[width=6.5cm,height=5.25cm]{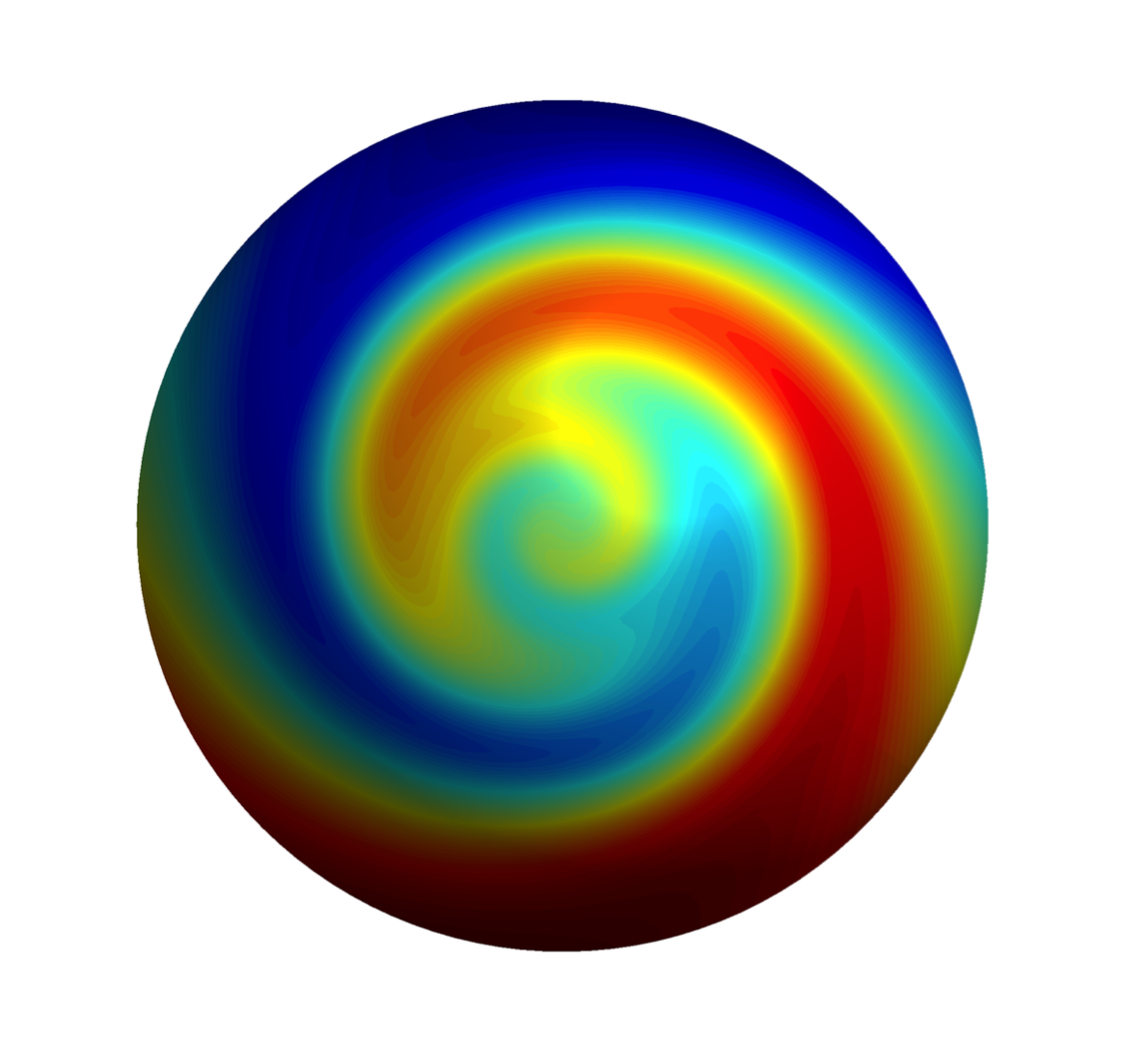}
    \includegraphics[width=6.5cm,height=5.25cm]{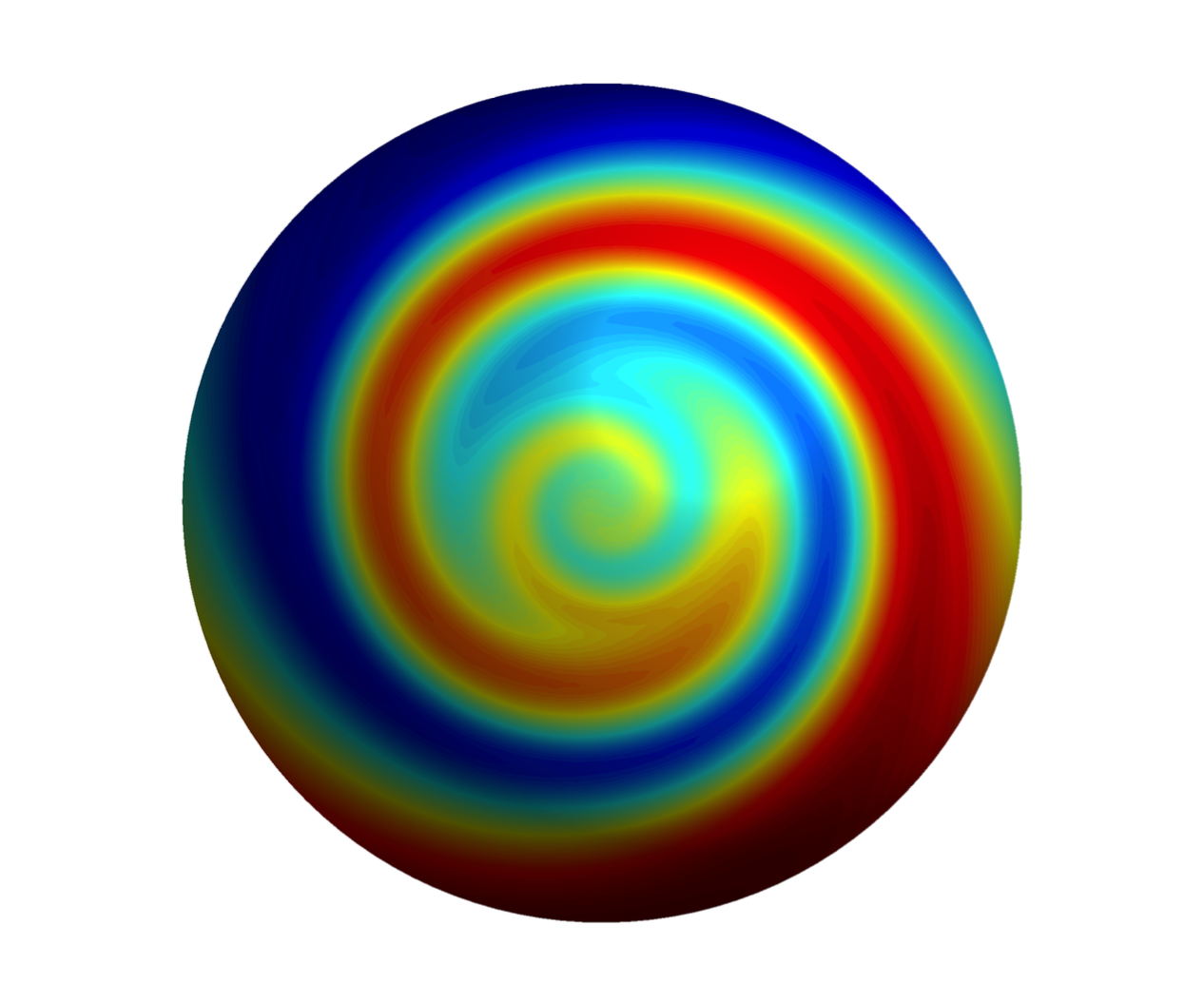}
    \caption{The 3D simulation of the vortex roll-up via MKLS approximation
        for the second test (given in Subsection \ref{62})  at $t=3$ (top left), $t=6$ (top right), and $t=9$ (bottom).}
    \label{fig-4}
\end{figure}

As the second standard test for transport equation on the unit
sphere, we have considered vortex roll-up test case, which is known
as deformational flow, and it models idealized cyclogenesis
\cite{fornberg2011stabilization,nair1999cascade}. We have solved the transport equation
(\ref{Eq-1}) with the following velocity field
\cite{fornberg2011stabilization,nair1999cascade}
$$v_{1}(\lambda,\theta)=\omega(\theta)\cos(\theta),\,\,\,\,\,\,v_{2}(\lambda,\theta)=0,$$
where
\[\omega \left( \theta  \right) = \left\{ \begin{array}{l}
\dfrac{{3\sqrt 3 }}{{2\rho \left( \theta  \right)}}\sec {h^2}
\left( {\rho \left( \theta  \right)} \right)\tanh \left( {\rho \left( \theta
\right)} \right),\,\,\,\,\,\,\,\,\,\,\rho \left( \theta  \right) \ne 0,\\
0,\,\,\,\,\,\,\,\,\,\,\,\,\,\,\,\,\,\,\,\,\,\,\,\,\,\,\,\,\,\,\,\,\,\,\,\,\,\,\,\,\,\,\,\,
\,\,\,\,\,\,\,\,\,\,\,\,\,\,\,\,\,\,\,\,\,\,\,\,\,\,\,\,\,\,\,\,\,\,\,\,\,\,\,\,\,\,\,\rho
\left( \theta  \right) \ne 0,\,
\end{array} \right.\]
where $\rho(\theta)=\rho_{0}\cos(\theta)$, and $\rho_{0}$ controls
the radial extent of the vortex \cite{fornberg2011stabilization,nair1999cascade}.
The analytical solution of this test is given as follows \cite{fornberg2011stabilization,nair1999cascade}
$$u(\lambda,\theta,t)=1-\tanh\left( \dfrac{\rho(\theta)}{\zeta}
\sin(\lambda-\omega(\theta)t)\right),\,\,\,\, t\geq 0.$$ For the
simulations reported here, we have chosen $\rho_{0}=3$ and
$\zeta=5$, which were considered previously in
\cite{fornberg2011stabilization}. All required parameters in two
approximations are chosen as the previous test, and $\Delta
t=1/1000$. In Figures \ref{fig-3} and \ref{fig-4}, we showed the
numerical solutions of Eq. (\ref{Eq-1}) on the unit sphere via
$N=19600$ PTS points at different time levels $t=3,6$ and $t=9$ by
GMLS and MKLS approximations. The results obtained in this test are
in agreement with reported results in
\cite{fornberg2011stabilization}. \DO{ In Table \ref{Table1-2}, the
used CPU time in both techniques for different values of $N$ during
simulations are given.} In Tables \ref{Table-3} and \ref{Table-4},
$\ell_{2}$ errors are reported at $T=3$ for different values of $N$
and $\Delta t=T/1000$ via GMLS and MKLS approximations,
respectively. \TODO{The results given here are in good agreement
with those reported in \cite{fornberg2011stabilization}, and we can
see almost the same accuracy using both approximations.}

\begin{table}
\begin{center}
\begin{tabular}{lllllllllllllllllll}
  \hline
  $\textbf{Method}$&&&&$N$ &&&&$\textbf{CPU time}\,(s)$     \\  
  \hline
\vspace{0.1cm} $\textbf{GMLS}$&&&&$1600$&&&& $28.50$   \\ \vspace{0.1cm}
 &&&&$6400$&&&& $108.86$ \\ \vspace{0.1cm}
 &&&&$19600$&&&& $256.68$   \\
 \hline
\vspace{0.1cm} $\textbf{MKLS}$&&&&$1600$&&&& $20.94$   \\ \vspace{0.1cm}
 &&&&$6400$&&&& $68.64$   \\ \vspace{0.1cm}
 &&&&$19600$&&&& $178.37$   \\
  \hline
\end{tabular}
\caption{\DO{The used CPU used time with different values $N$ \\ for the BiCGSTAB method
for the second test problem.}}\label{Table1-2}
\end{center}
\end{table}

\begin{table}
\begin{center}
\begin{tabular}{llllllllllllll}
  \hline
  &&\multicolumn{2}{l}{PTS}&&\multicolumn{2}{l}{ME} \\
  \cline{3-4} \cline{6-7}
  $N$&&$\ell_{2}$ &&&& $\ell_{2}$ \\
  \hline
\vspace{0.1cm} $400$  && $2.25\e-2$     &&&& $2.09\e-2$      \\ \vspace{0.1cm}
 $1600$  && $3.51\e-3$    &&&& $3.46\e-3$    \\ \vspace{0.1cm}
 $6400$ && $5.22\e-4$   &&&& $7.78\e-4$    \\ \vspace{0.1cm}
 $16641$ && $1.97\e-4$  &&&& $1.92\e-4$  \\
  \hline
\end{tabular}
\caption{The $\ell_{2}-$error for different values $N$ \\ in GMLS approximation for the second test.}\label{Table-3}
\end{center}
\end{table}

\begin{table}
\begin{center}
\begin{tabular}{llllllllllllll}
  \hline
  &&\multicolumn{2}{l}{PTS}&&\multicolumn{2}{l}{ME} \\
  \cline{3-4} \cline{6-7}
  $N$&&$\ell_{2}$ &&&& $\ell_{2}$ \\
  \hline
  \vspace{0.1cm}$400$  && $4.05\e-2$     &&&& $4.17\e-2$      \\ \vspace{0.1cm}
 $1600$  && $1.41\e-2$    &&&& $1.31\e-2$    \\ \vspace{0.1cm}
 $6400$ && $3.59\e-3$   &&&& $1.75\e-3$    \\ \vspace{0.1cm}
 $16641$ && $7.48\e-4$  &&&& $2.10\e-4$  \\
  \hline
\end{tabular}
\caption{The $\ell_{2}-$error for different values $N$ \\
in MKLS approximation for the second test.}\label{Table-4}
\end{center}
\end{table}

\subsection{Deformational flow test}\label{63}
In this part, we consider the following test, which is
known as deformational flow \cite{nair2010class}. We have considered the transport equation (\ref{Eq-1})
with the following velocity field
$$v_{1}(\lambda,\theta,t)=2\sin^2(\lambda)\sin(2\theta)\cos(\pi t/T),\,\,\,\,\,\,
v_{2}(\lambda,\theta,t)=2\sin(2\lambda)\cos(\theta)\cos(\pi t/T),$$
which is
non-divergent flow \cite{nair2010class}. The initial condition for this test is given as follows \cite{nair2010class}
\[u(\lambda ,\theta,t=0 ) = \left\{ \begin{array}{l}
0.1 + 0.9{u_1}(\lambda ,\theta ),\,\,\,\,\,\,\,\,\,\,\,{r_1}(\lambda ,\theta ) < r,\\
0.1 + 0.9{u_2}(\lambda ,\theta ),\,\,\,\,\,\,\,\,\,\,\,{r_2}(\lambda ,\theta ) < r,\\
0.1,\,\,\,\,\,\,\,\,\,\,\,\,\,\,\,\,\,\,\,\,\,\,\,\,\,\,\,\,\,\,\,\,\,\,\,\,\,\,\,\,\,\,\,\,\,\,o.w,
\end{array} \right.\]
where
\[{u_1}(\lambda ,\theta ) = \dfrac{1}{2}\left( {1 + \cos \left( {\frac{{\pi {r_1}(\lambda ,\theta )}}{r}} \right)}
\right),\,\,\,\,\,\,\,\,\,\,\,\,\,{u_2}(\lambda ,\theta ) =
 \dfrac{1}{2}\left( {1 + \cos \left( {\frac{{\pi {r_2}(\lambda ,\theta )}}{r}} \right)} \right),\,\,\,\]
and
\[\begin{array}{l}
{r_1}(\lambda ,\theta ) = \arccos \left( {\sin ({\theta _1})\,\sin (\theta ) + \cos ({\theta _1})\,\cos (\theta )\cos (\lambda  - {\lambda _1})} \right)\,,\\\\
{r_2}(\lambda ,\theta ) = \arccos \left( {\sin ({\theta _2})\,\sin (\theta ) + \cos ({\theta _2})\,\cos (\theta )
\cos (\lambda  - {\lambda _2})} \right).
\end{array}\]

\begin{figure}[t!]
    \centering
    \includegraphics[width=6.5cm,height=5.25cm]{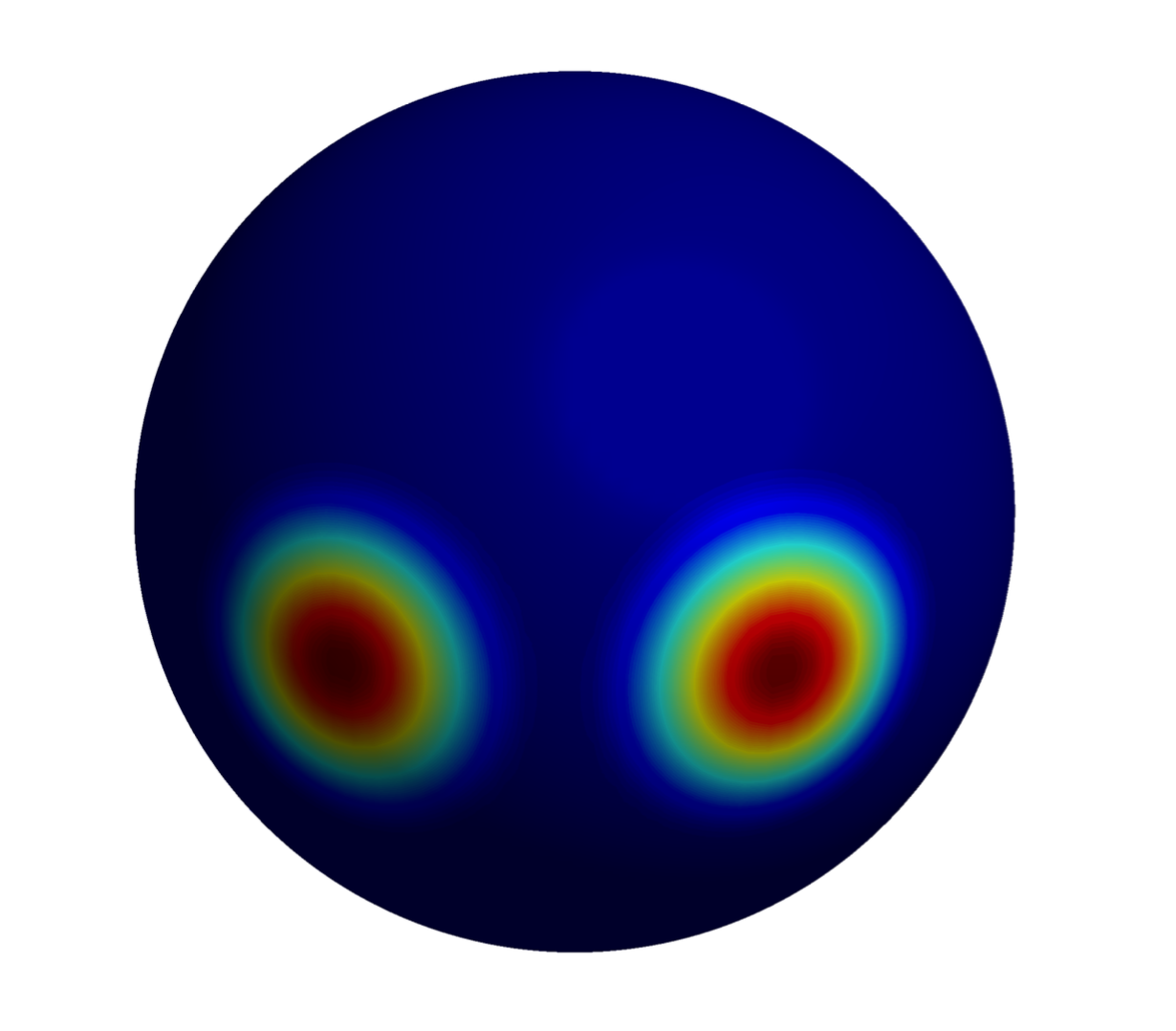}
    \caption{The initial condition of the deformational flow in  the third test (given in Subsection \ref{63}).}
    \label{fig-4-1}
\end{figure}

\begin{figure}[t!]
    \centering
    \includegraphics[width=6.5cm,height=5.25cm]{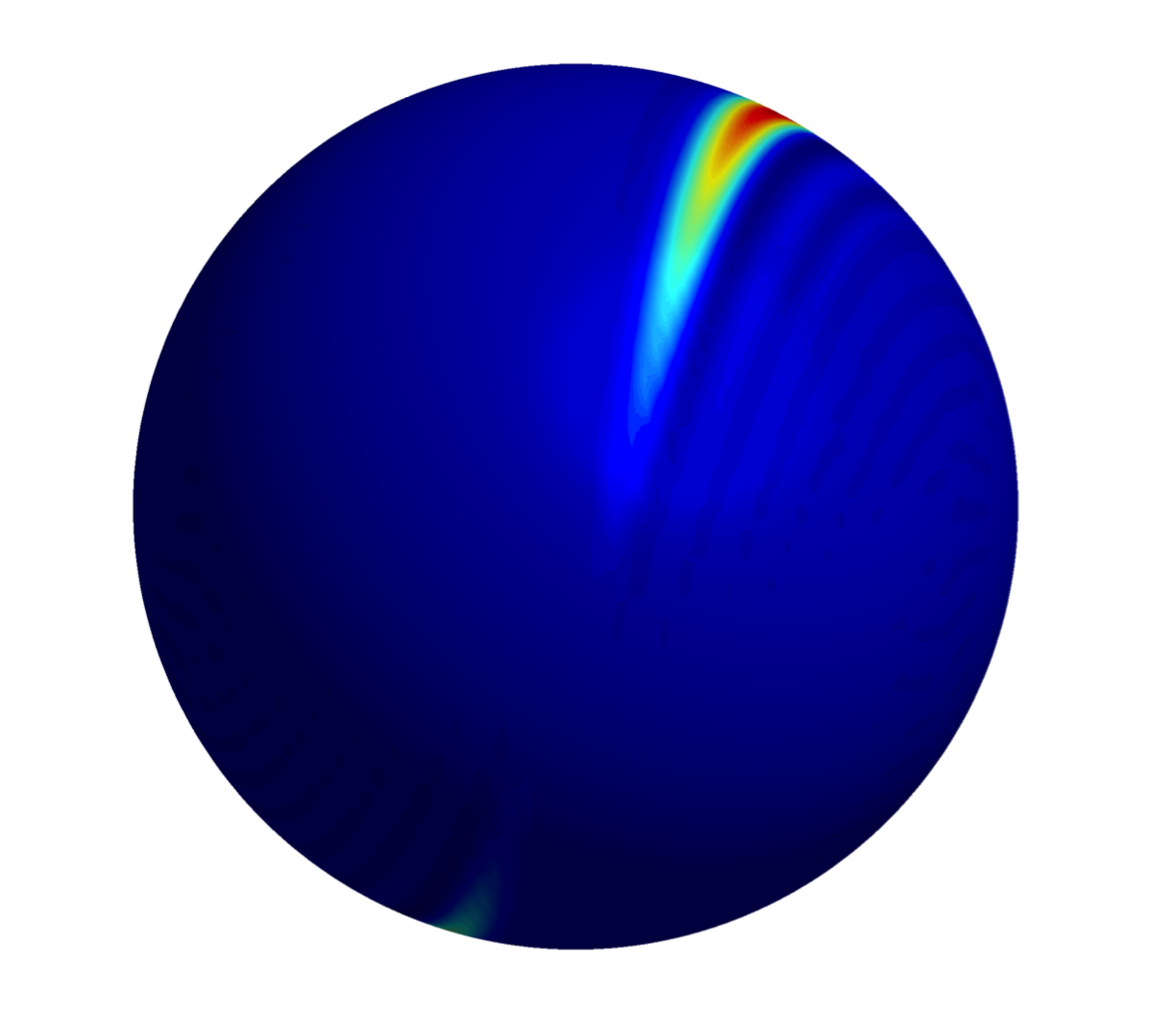}\hspace{1cm}
    \includegraphics[width=6.5cm,height=5.25cm]{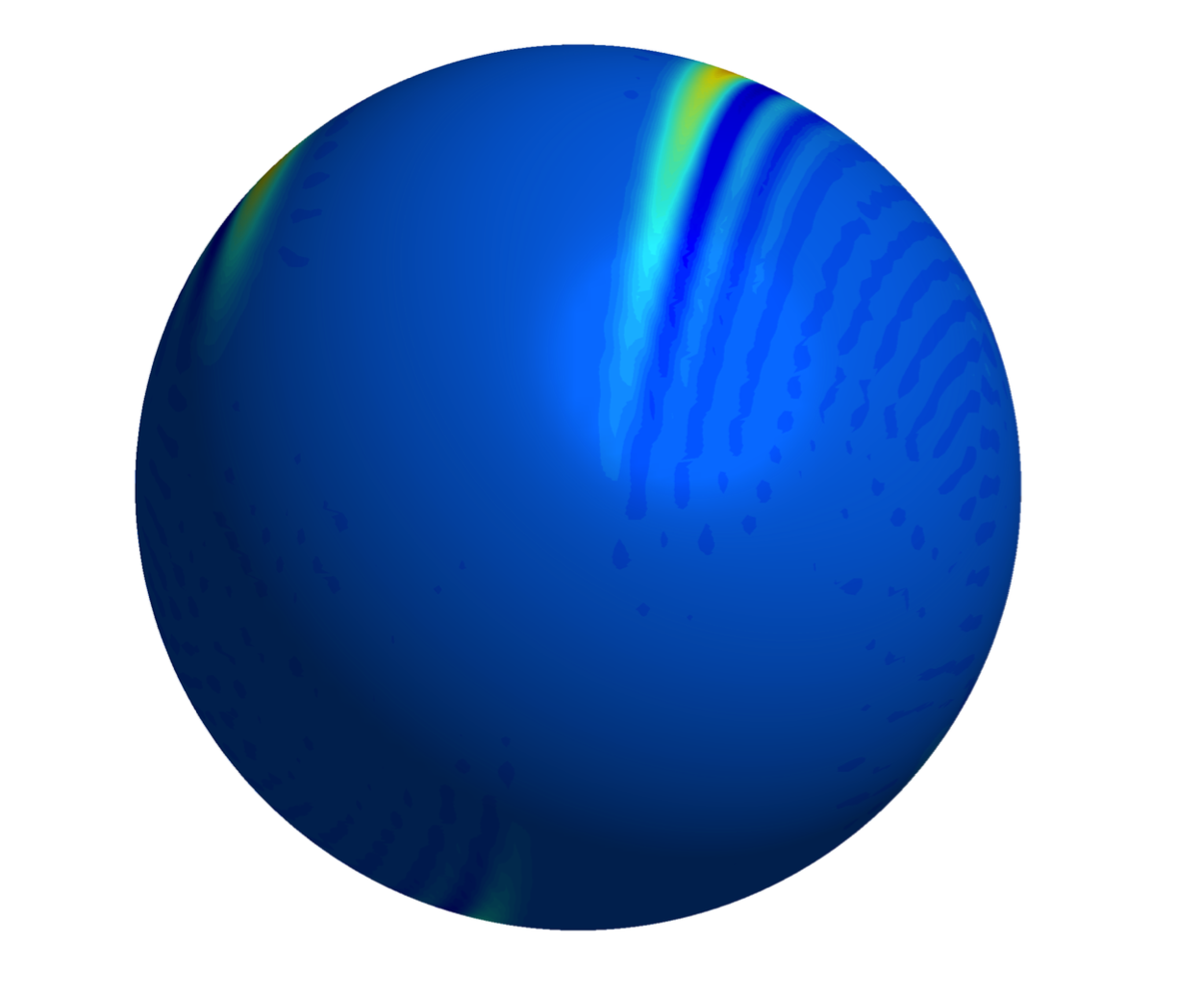}\vspace{1cm}
    \includegraphics[width=6.5cm,height=5.25cm]{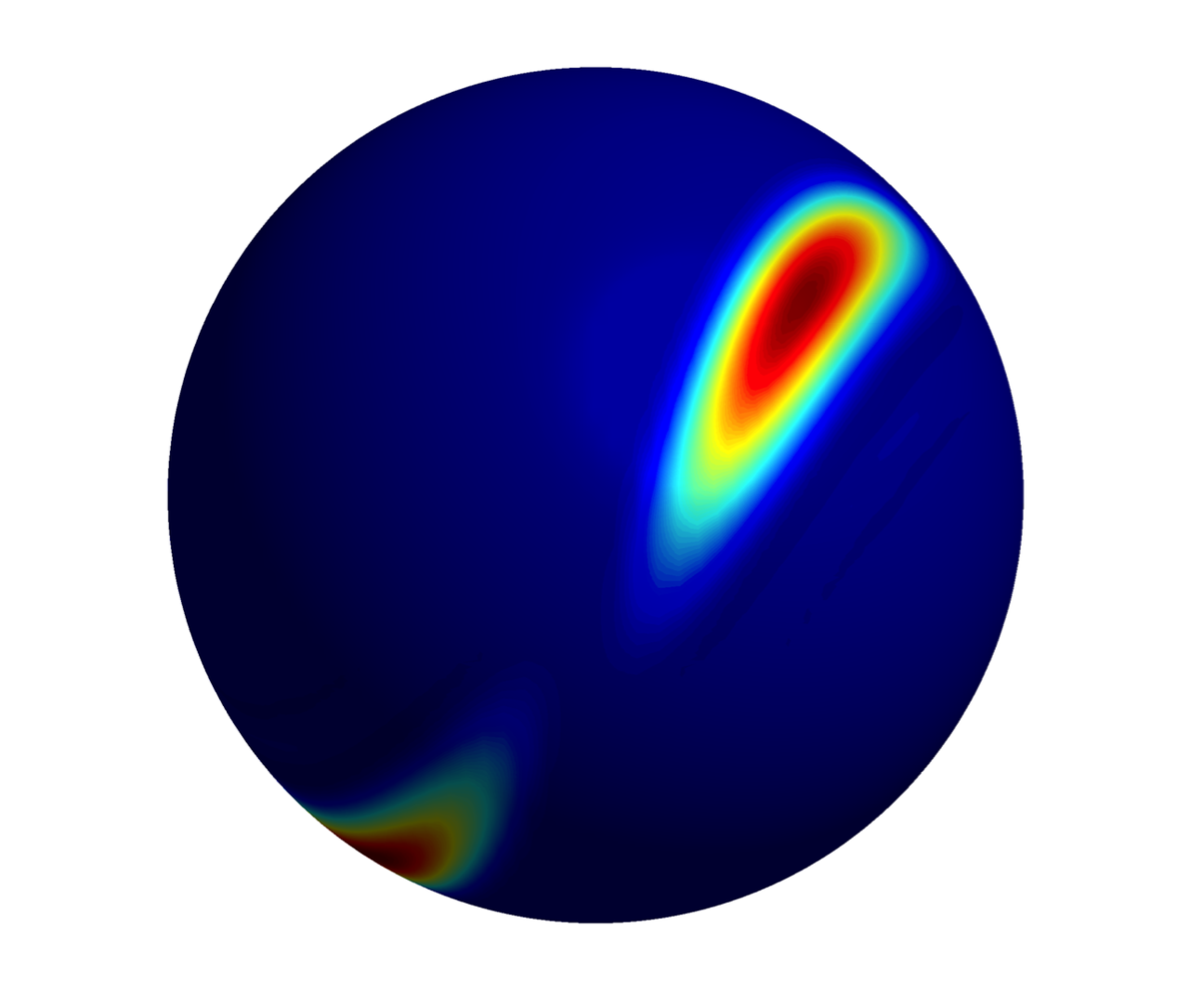}\hspace{1cm}
    \includegraphics[width=6.5cm,height=5.25cm]{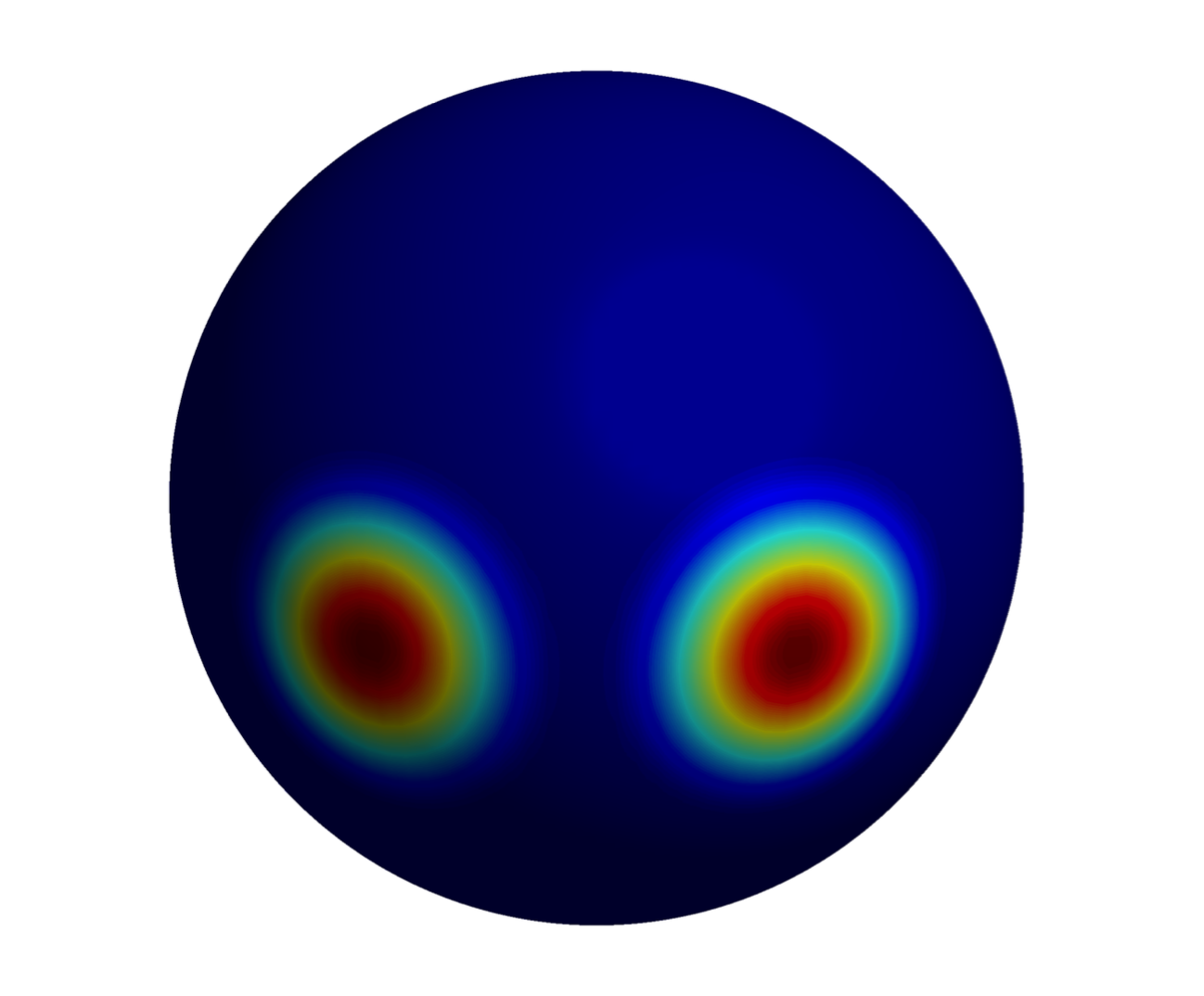}
    \caption{The 3D simulation of the deformational flow via GMLS approximation
        for the third test (given in Subsection \ref{63}) at $t=T/4$ (top left),  $t=T/2$ (top right),  $t=9T/2$ (bottom left), and  $t=T$ (bottom right).}
    \label{fig-5}
\end{figure}

\begin{figure}[t!]
    \centering
    \includegraphics[width=6.5cm,height=5.25cm]{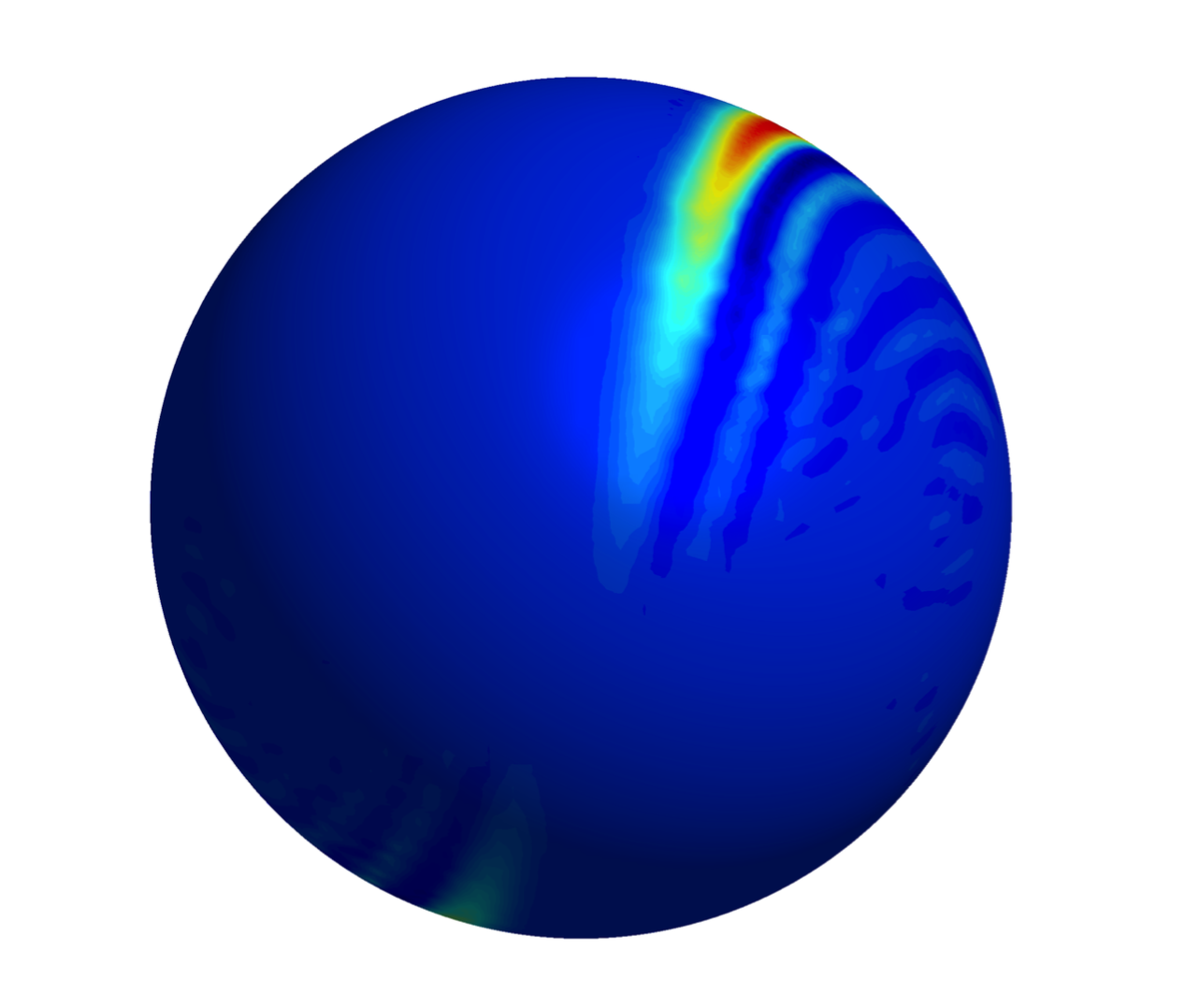}\hspace{1cm}
     \includegraphics[width=6.5cm,height=5.25cm]{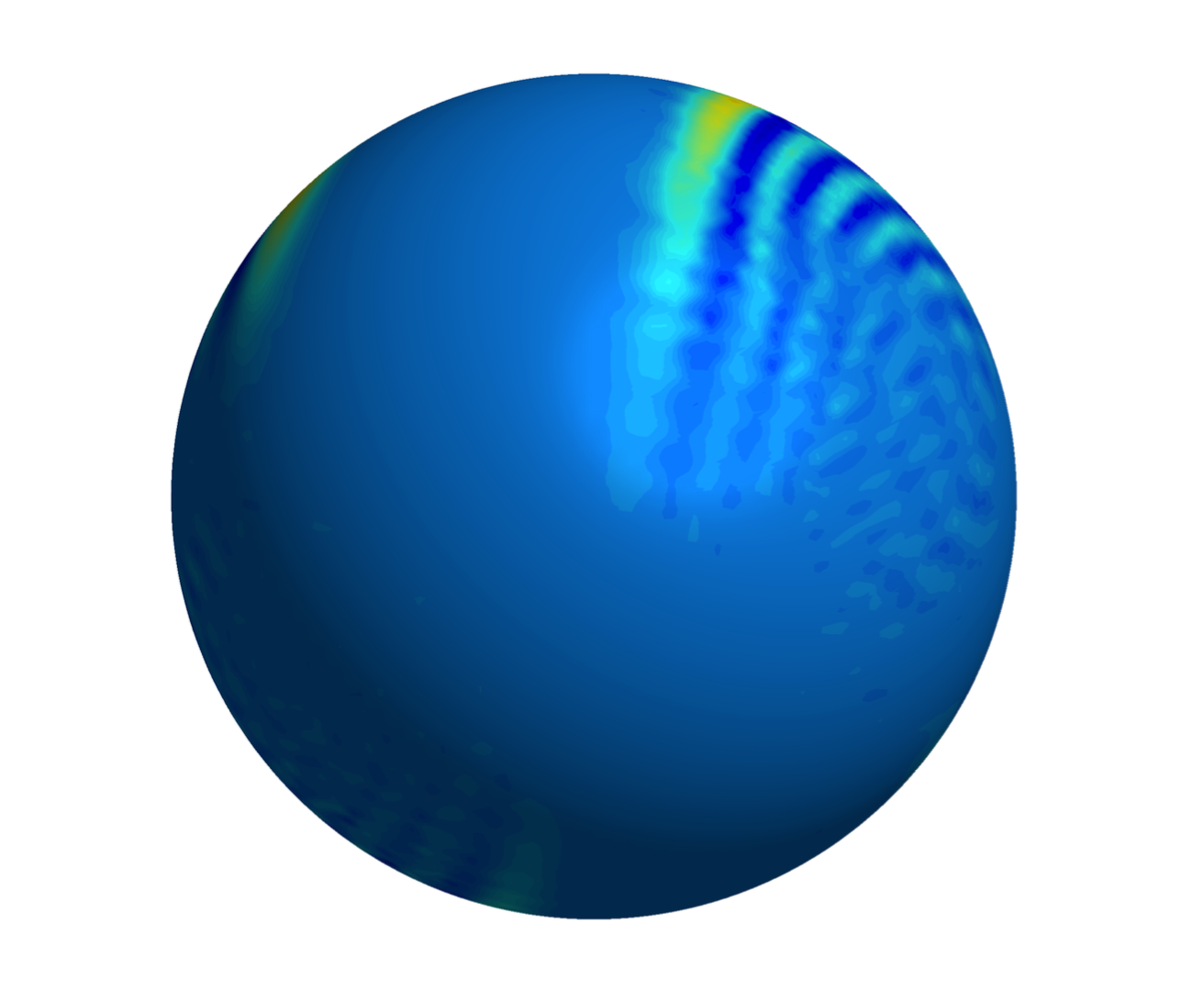}\vspace{1cm}
      \includegraphics[width=6.5cm,height=5.25cm]{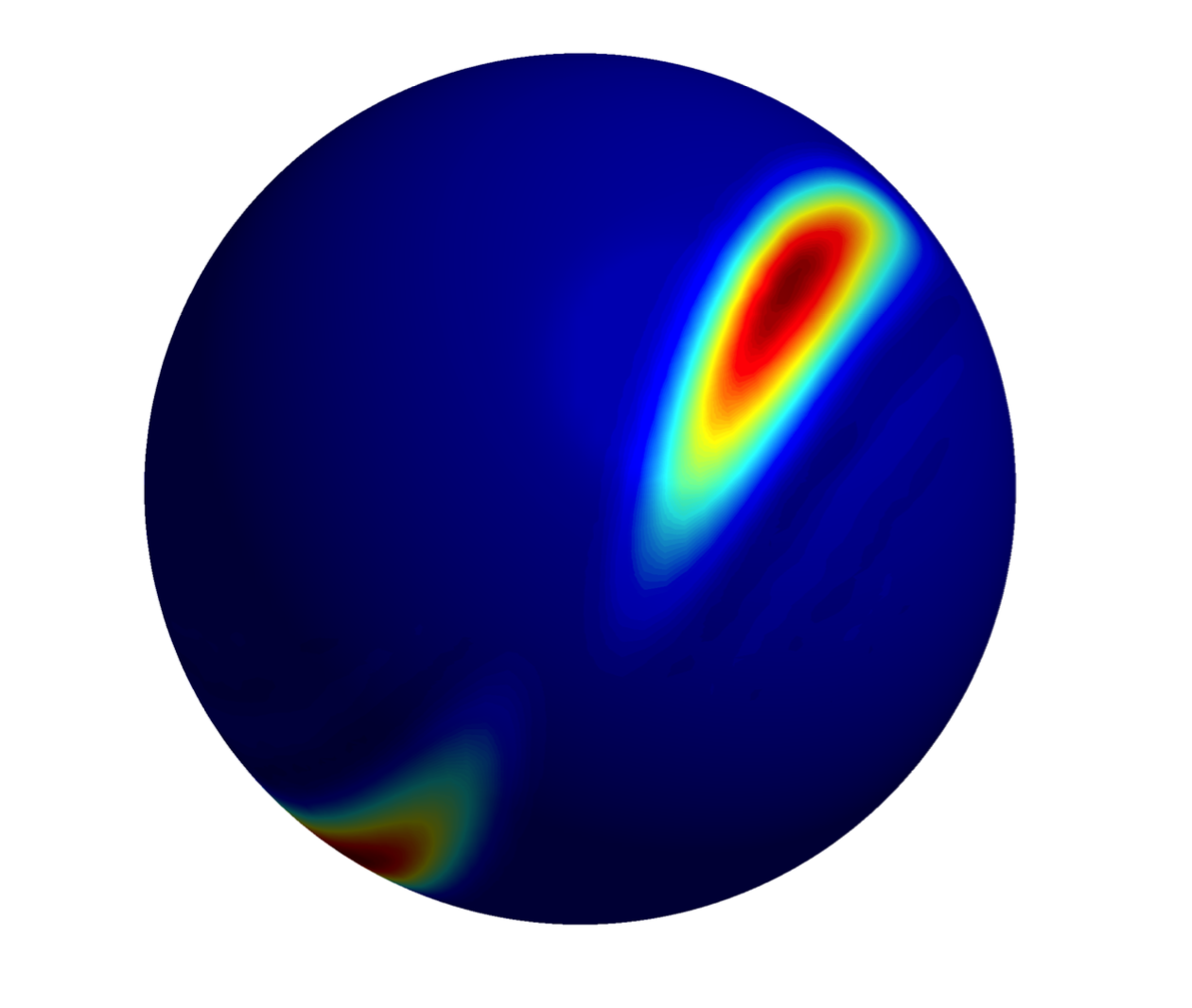}\hspace{1cm}
       \includegraphics[width=6.5cm,height=5.25cm]{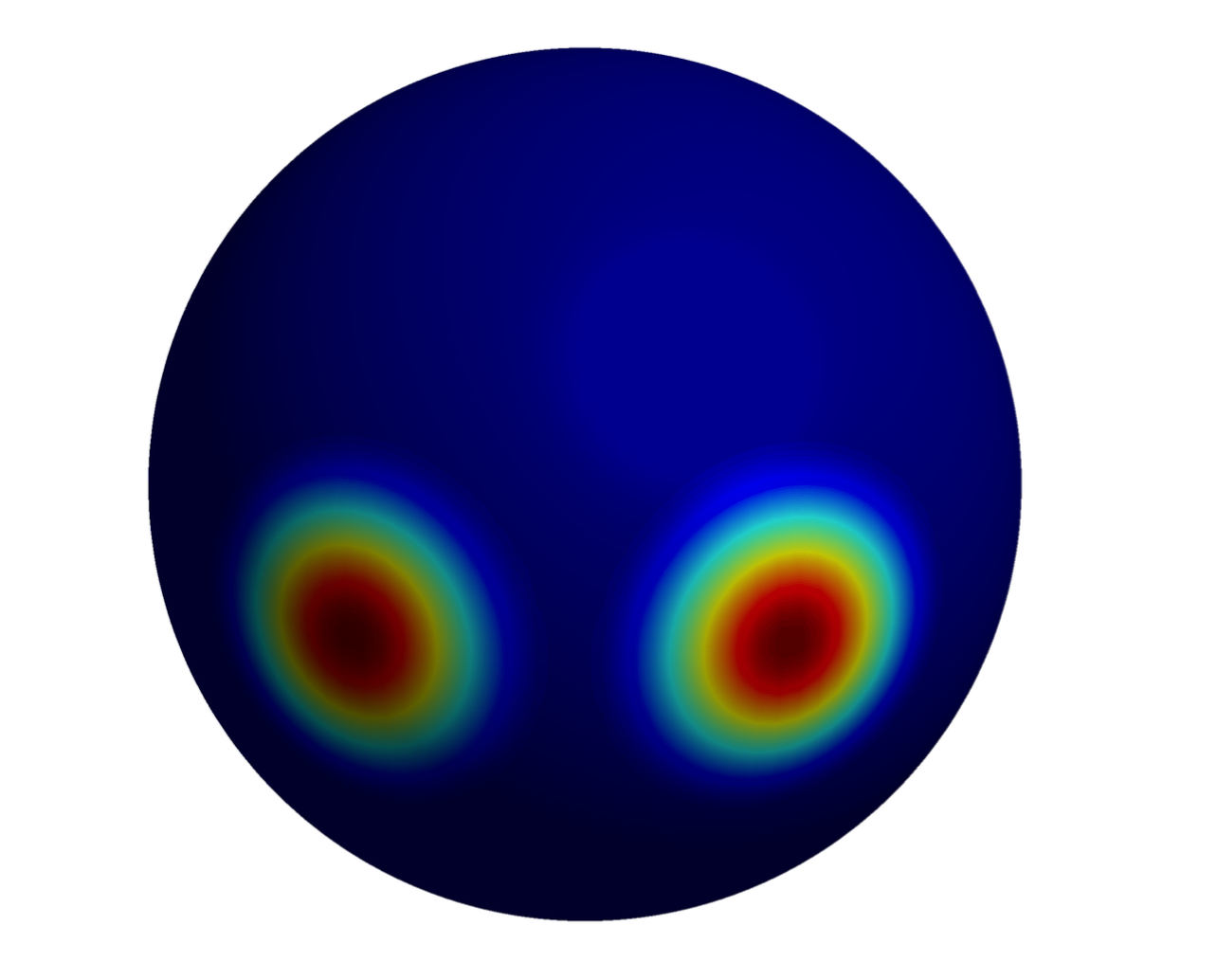}
    \caption{The 3D simulation of the deformational flow via MKLS approximation
        for the third test (given in Subsection \ref{63}) at $t=T/4$ (top left),  $t=T/2$ (top right),  $t=9T/2$ (bottom left), and  $t=T$ (bottom right).}
    \label{fig-6}
\end{figure}

In the above formulations, $(\lambda_{1},\theta_{1})=(5\pi/6,0)$ and
$(\lambda_{2},\theta_{2})=(7\pi/6,0)$ are the centers of two cosine
bells \cite{nair2010class}. This test shows that the flow field will
be deformed at $t=2.5$, and then with return to its initial position
(see Figure \ref{fig-4-1}) at $T=5$ \cite{nair2010class}. In Figure
\ref{fig-5}, the numerical solution of $u$ at different time levels
$t=T/4,T/2,9T/2$ and $t=T$ via GMLS approximation, where $T=5$ and
$\Delta t=1/400$ and $N=6400$ ME points is shown. Furthermore, the
same simulations are obtained in Figure \ref{fig-6} with MKLS
approximation.

\TODO{As observed and expected in these figures, two cosine bells
considered at the initial (Figure \ref{fig-4-1}) are deformed at
$t=T/2$, and they are returned to their initial positions at $t=T$.}
\DO{ In Table \ref{Table1-3}, the used CPU time in GMLS and MKLS
approximations for different values of $N$ during simulations are
shown.} In Table \ref{Table-5}, $\ell_{2}$ errors are given via GMLS
approximation using PTS and ME points on the
 unit sphere. In Table \ref{Table-6}, $\ell_{2}$ errors
 obtained from the implementation of
  MKLS technique
with set of points considered are reported and \TO{different values
time step such that for $N=400$, $\Delta t=1/100$, for $N=1600$,
$\Delta t=1/200$, for $N=6400$, $\Delta t=1/400$, for $N=16641$,
$\Delta t=1/800$.} \TODO{Also in these tables, the accuracy in both
approximations is almost same.}

\begin{table}
\begin{center}
\begin{tabular}{lllllllllllllllllll}
  \hline
  $\textbf{Method}$&&&&$N$ &&&&$\textbf{CPU time}\,(s)$     \\
  \hline
  \vspace{0.1cm}$\textbf{GMLS}$&&&&$1600$&&&& $5.98$   \\ \vspace{0.1cm}
 &&&&$6400$&&&& $19.28$   \\ \vspace{0.1cm}
 &&&&$10000$&&&& $28.11$ \\
 \hline
\vspace{0.1cm} $\textbf{MKLS}$&&&&$1600$&&&& $5.58$   \\ \vspace{0.1cm}
 &&&&$6400$&&&& $19.25$   \\ \vspace{0.1cm}
 &&&&$10000$&&&& $26.67$   \\
  \hline
\end{tabular}
\caption{\DO{The used CPU time   with different values $N$\\ for the BiCGSTAB method
for the third test.}}\label{Table1-3}
\end{center}
\end{table}

\begin{table}
\begin{center}
\begin{tabular}{llllllllllllll}
  \hline
  &&\multicolumn{2}{l}{PTS}&&\multicolumn{2}{l}{ME} \\
  \cline{3-4} \cline{6-7}
  $N$&&$\ell_{2}$ &&&& $\ell_{2}$ \\
  \hline
\vspace{0.1cm} $400$  && $3.34\e-3$     &&&& $3.43\e-3$      \\ \vspace{0.1cm}
 $1600$  && $1.11\e-3$    &&&& $1.15\e-3$    \\ \vspace{0.1cm}
 $6400$ && $2.76\e-4$   &&&& $3.96\e-4$    \\ \vspace{0.1cm}
 $16641$ && $8.36\e-5$  &&&& $1.11\e-4$  \\
  \hline
\end{tabular}
\caption{The $\ell_{2}-$error for different values $N$ \\ in GMLS
approximation for the third test.}\label{Table-5}
\end{center}
\end{table}

\begin{table}
\begin{center}
\begin{tabular}{llllllllllllll}
  \hline
  &&\multicolumn{2}{l}{PTS}&&\multicolumn{2}{l}{ME} \\
  \cline{3-4} \cline{6-7}
  $N$&&$\ell_{2}$ &&&& $\ell_{2}$ \\
  \hline
\vspace{0.1cm} $400$  && $1.45\e-3$     &&&& $1.46\e-3$      \\ \vspace{0.1cm}
 $1600$  && $8.75\e-4$    &&&& $8.81\e-4$    \\ \vspace{0.1cm}
 $6400$ && $2.53\e-4$   &&&& $2.55\e-4$    \\ \vspace{0.1cm}
 $16641$ && $2.00\e-4$  &&&& $1.51\e-4$  \\
  \hline
\end{tabular}
\caption{The $\ell_{2}-$error for different values $N$ \\ in MKLS
approximation for the third test.}\label{Table-6}
\end{center}
\end{table}

\TO{\subsection{Comparison between two proposed approximations with
other methods}
 In this section, we compare GMLS and MKLS approximations with other numerical
 methods, which were applied to solve the transport equation on the
 sphere in the literature, we have considered CSLAM method
 \cite{lauritzen2012standard},
 DG method \cite{nair2010class}, RBF-FD
 technique \cite{fornberg2011stabilization}, local and global RBF
 approaches \cite{shankar2018mesh} and RBF-PU method
 \cite{shankar2018mesh}. We have used the deformational flow test
 for the cosine bell, which is given in Subsection \ref{63}. A
 comparison between all mentioned methods is done in Table \ref{Table-7} due to degrees of freedom (DOF)
 (the unknown coefficients
 related to each method), and
 time steps ($\Delta t$) by computing relative $\ell_{2}$ error. The errors
 reported here for GMLS and MKLS approximations are computing via PTS points on $\mathbb{S}^2$ due to
 formula (\ref{Integration-1}), which approximates $\ell_{2}$
 error. It also should be noted that,
 in \cite{shankar2018mesh}, for computing the surface integral, the sixth--order
 kernel--based meshfree quadrature method has been used \cite{fuselier2014kernel}.
 Besides, the results reported here for
 CSLAM, DG, RBF-FD, local RBF, global RBF, and RBF-PU methods
 are taken from \cite[Table 2, Subsection 4.4]{shankar2018mesh}. }

\begin{table}
\begin{center}
\begin{tabular}{llllllllllllllllllllllllllll}
  \hline
  \cline{3-4} \cline{6-7}
  \textbf{Method}&&&$\Delta t$ &&&& $\textbf{DOF}$&&&& \textbf{Relative} $\ell_{2}$ \textbf{error} \\
  \hline
 \vspace{0.1cm} \textbf{GMLS}  &&& $5/2400$     &&&& $6400$ &&&& $2.84\e-4$      \\ \vspace{0.1cm}
 \textbf{MKLS}  &&& $5/2400$    &&&& $6400$ &&&& $2.81\e-4$    \\ \vspace{0.1cm} 
 \textbf{CSLAM} \cite{lauritzen2012standard}  &&& $5/240$    &&&& $86400$ &&&& $6.00\e-3$    \\ \vspace{0.1cm}
 \textbf{DG}, $p=3$ \cite{nair2010class} &&& $5/2400$   &&&& $38400$ &&&& $1.39\e-2$    \\ \vspace{0.1cm}
 \textbf{RBF-FD}, $n=84$ \cite{nair2010class} &&& $5/900$  &&&& $23042$ &&&& $1.17\e-2$  \\ \vspace{0.1cm}
 \textbf{Local RBF}, $n=84$ \cite{shankar2018mesh} &&& $5/35$  &&&& $23042$ &&&& $3.45\e-3$  \\ \vspace{0.1cm}
 \textbf{RBF-PU}, $n=84$ \cite{shankar2018mesh}&&& $5/35$  &&&& $23042$ &&&& $3.63\e-3$  \\ \vspace{0.1cm}
 \textbf{Global RBF} \cite{shankar2018mesh}&&& $5/45$  &&&& $15129$ &&&& $5.10\e-3$  \\
  \hline
\end{tabular}
\caption{\TO{A comparison between the presented approximations and
other numerical methods for deformational flow test (cosine bell)
due to DOF and time steps. The compared methods with GMLS and MKLS
techniques are CSLAM \cite{lauritzen2012standard}, which is based on
a cubed sphere grid. DG is the discontinuous Galerkin scheme
\cite{nair2010class} with $p=3$ degree polynomials (fourth--order
accurate) and the cubed sphere grid. RBF-FD is the mesh-free
Eulerian scheme \cite{fornberg2011stabilization} via $n=84$ points
in each local domain. Semi-Lagrangian Local RBF, RBF-PU and global
RBF methods, which are applied in
\cite{shankar2018mesh}.}}\label{Table-7}
\end{center}
\end{table}

\newpage


\TO{\section{Concluding remarks}\label{Sec-7} In this paper, two new
techniques, i.e., GMLS and MKLS have been applied to approximate the
spatial variables of a transport equation on the sphere in spherical
coordinates. The time variable of the model is discretized by a
second-order backward differential formula. The resulting of 
fully discrete scheme is a linear system of algebraic equations per
time step, which is solved efficiently by a BiCGSTAB method with
zero-fill ILU preconditioner. To ensure the ability of the proposed approaches, we have solved three important test cases namely solid body rotation, vortex roll-up, and deformational flow, which are important examples in the numerical climate modeling community.  Both developed
techniques do not depend on a background mesh or
triangulation, which yields an easy implementation in solving of the
transport equation on the sphere.  Due to this feature, we have
obtained the numerical results using two different distribution
points, i.e., PTS and ME on the sphere. Furthermore, pole
singularities appeared in this equation have been omitted in
differentiation matrices due to the applied approximations.
 As formulated in our Algorithms
\ref{Algorithm1} and \ref{Algorithm2}, the procedure of
implementation in both approximations consists of two main parts,
which are the construction of required matrices due to Eqs.
(\ref{full-5}) and (\ref{full-6}) (Algorithm \ref{Algorithm1}) and
solving the obtained full-discrete scheme (Algorithm
\ref{Algorithm2}). As mentioned before, all implementations are done in MATLAB software by writing routines due to the presented
algorithms. The results and simulations reported here show that both
methods have the same accuracy, but the MKLS approximation depends on a
constant parameter, which should be controlled experimentally.
Besides, as shown in the first test, the GMLS approximation uses less
CPU time than MKLS approximation for constructing all matrices in
Algorithm \ref{Algorithm1}.
 We also have compared GMLS and MKLS approximations
with other methods in the literature, i.e., CSLAM, DG, RBF-FD, local
RBF, global RBF, and RBF-PU for deformational flow test due to time
steps and DOF. From these comparisons, we can observe that the
accuracy of GMLS and MKLS approximations with less number of points
(DOF) is better than other methods. According to the results and
discussions in this paper, the GMLS and MKLS approaches can be applied
easily to solve mathematical models in spherical geometries.}

%

\bibliographystyle{elsarticle-num}
\bibliography{MLMC}
\newpage

\end{document}